\documentclass[12pt]{amsart}
\usepackage{amssymb}
\usepackage{amsmath}
\usepackage{amsfonts}
\usepackage{mathrsfs}
\usepackage{graphicx}
\usepackage{color}
\usepackage[onehalfspacing]{setspace}
\usepackage{caption}
\usepackage{enumerate}
\usepackage[utf8]{inputenc}
\usepackage[charter,cal=cmcal]{mathdesign}
\usepackage{epigraph}
\numberwithin{equation}{section}
\theoremstyle{plain}

\begin{document}
\title[Local Asymptotic Minimax Estimation of Nonregular Parameters]{Local Asymptotic Minimax
Estimation of Nonregular Parameters with Translation-Scale Equivariant Maps}
\author[K. Song]{Kyungchul Song}

\address{Vancouver School of Economics, University of British Columbia, 6000 Iona Drive, Vancouver, BC, V6S 1L4, Canada}
\email{kysong@mail.ubc.ca}
\thanks{\today \\  This version is an update of my paper published in Journal of Multivariate Analysis (Song (2014)). The update is made minimally, only filling two gaps in the mathematical proofs. At the end of this paper, I attach a corrigendum that shows where the gaps are in the published paper. I thank Keisuke Hirano, Marcelo Moreira, Ulrich M\"{u}ller and Frank Schorfheide for valuable comments for the earlier version of the paper. I also thank Yoshiyasu Rai and Zheng Fang for pointing out the gaps in the mathematical proofs of the paper. This research was supported by
	the Social Sciences and Humanities Research Council of Canada.}
	
\begin{abstract}
{\footnotesize When a parameter of interest is defined to be a
nondifferentiable transform of a regular parameter, the parameter does not
have an influence function, rendering the existing theory of semiparametric
efficient estimation inapplicable. However, when the nondifferentiable
transform is a known composite map of a continuous piecewise linear map with
a single kink point and a translation-scale equivariant map, this paper
demonstrates that it is possible to define a notion of asymptotic optimality
of an estimator as an extension of the classical local asymptotic minimax
estimation. This paper establishes a local asymptotic risk bound and
proposes a general method to construct a local asymptotic minimax decision.}

{\footnotesize \ }

{\footnotesize \noindent \textsc{Key words.} Nonregular Parameters;
Translation-Scale Equivariant Transforms; Semiparametric Efficiency; Local
Asymptotic Minimax Estimation. \newline
}

{\footnotesize \noindent \textsc{AMS Classification. 62C05, 62C20.}}
\end{abstract}

\maketitle

\section{Introduction}

This paper investigates the problem of optimal estimation of a parameter $%
\theta \in \mathbf{R}$ which takes the following form:%
\begin{equation}
\theta =(f\circ g)(\mathbf{\beta }),  \label{th}
\end{equation}%
where $\mathbf{\beta }\in \mathbf{R}^{d}$ is a regular parameter for which a
semiparametric efficiency bound is well defined, $g$ is a translation-scale
equivariant map, and $f$ is a continuous piecewise linear map with a single
kink (i.e., nondifferentiability) point.

Examples abound, including $\max \{\beta _{1},\beta _{2},\beta _{3}\}$, $%
\max \{\beta _{1},0\}$, $|\beta _{1}|$, $|\max \{\beta _{1},\beta _{2}\}|$,
etc., where $\mathbf{\beta }=(\beta _{1},\beta _{2},\beta _{3})$ is a
regular parameter, i.e., a parameter which is differentiable in the
underlying probability. Applications where such parameters arise are
numerous. We give two specific examples.\bigskip

\noindent \textsc{Example 1 (Maximal Average Treatment Effects):} Suppose
that $X$ is an observed discrete covariate and $D\in \{0,1,2,\cdot \cdot
\cdot ,J\}$ is a treatment indicator, where $D=j$ for $j>0$ indicates
treatment by method $j$, and $D=0$ indicates no treatment. Let us assume
that the vector of potential outcomes $(Y_{0},Y_{1},\cdot \cdot \cdot
,Y_{J}) $ are conditionally independent from $X$ given $D$, and that $%
P\{D=j|X=x\}\in (0,1)$ for all $j=0,1,2,\cdot \cdot \cdot ,J$ and $x$ in the
support of $X$. The researcher observes $(Y,D,X)$, where $%
Y=\sum_{j=0}^{J}Y_{j}1\{D=j\}$, but does not observe $(Y_{j})_{j=0}^{J}$.
Then the average treatment effect for method $j$ for group with $X=x$ is
identified by%
\begin{equation*}
\beta _{j}=\mathbf{E}[Y|X=x,D=j]-\mathbf{E}[Y|X=x,D=0].
\end{equation*}%
One of the examples considered by Hirano and Porter (2012) was 
\begin{equation*}
\theta =\max_{1\leq j\leq J}\beta _{j}\text{,}
\end{equation*}%
that is, the maximum treatment effect that is possible using the $J$
methods. $\blacksquare $\bigskip

\noindent \textsc{Example 2\ (Bounds for Treatment Effects under
Monotonicity):} Let $Y_{j}$ be the potential outcome variables taking values
from $[K_{0},K_{1}]$ with known constants, $K_{0}$ and $K_{1}$, and $D$ a
treatment indicator as in Example 1. Suppose that $X$ is an observed
discrete random variable taking values in $\{x_{1},\cdot \cdot \cdot
,x_{M}\} $, $x_{1}\leq x_{2}\leq \cdot \cdot \cdot \leq x_{M-1}\leq x_{M}$,
such that $\mathbf{E}[Y_{j}|X=x]\geq \mathbf{E}[Y_{j}|X=x^{\prime }]$
whenever $x\geq x^{\prime }$ for all $j=0,1,\cdot \cdot \cdot ,J$. The
parameter of interest is the conditional outcome $\mathbf{E}[Y_{j}|X=x]$ for
treatment method $j$. The researcher observes $(Y,D,X)$ with $%
Y=\sum_{j=0}^{J}Y_{j}1\{D=j\}$ as before. Manski and Pepper (2000) showed
that in this set-up, the conditional outcome is interval identified as
follows:%
\begin{equation*}
\max_{1\leq k\leq m}\beta _{j,k}(K_{0})\leq \mathbf{E}[Y_{j}|X=x_{m}]\leq
\min_{m\leq k\leq M}\beta _{j,k}(K_{1}),
\end{equation*}%
where 
\begin{equation*}
\beta _{j,k}(K)=\mathbf{E}\left[ Y|X=x_{k},D=j\right] P\left\{
D=j|X=x_{k}\right\} +K\cdot P\left\{ D\neq j|X=x_{k}\right\} \text{.}
\end{equation*}%
Then the upper bound parameter $\theta _{U}=\min_{m\leq k\leq M}\beta
_{j,k}(K_{1})$ and the lower bound parameter $\theta _{L}=\max_{1\leq k\leq
m}\beta _{j,k}(K_{0})$ are examples of $\theta $ in (\ref{th}). Such a bound
frequently arises in economics literature (e.g. Haile and Tamer (2003) for
bidders' valuations in English auctions.) $\blacksquare $\bigskip

In contrast to the ease with which a parameter of the form in (\ref{th})
arises in applied researches, a formal analysis of the optimal estimation
problem has remained a challenging task. One might consistently estimate $%
\theta $ by using plug-in estimator $\hat{\theta}=f(g(\mathbf{\hat{\beta}}))$%
, where $\mathbf{\hat{\beta}}$ is a $\sqrt{n}$-consistent estimator of $%
\mathbf{\beta }$. However, there have been concerns about the asymptotic
bias that such an estimator carries, and some researchers have proposed ways
to reduce the bias (Manski and Pepper (2000), Haile and Tamer (2003),
Chernozhukov, Lee, and Rosen (2013)). However, Doss and Sethuraman (1989)
showed that a sequence of estimators of a parameter for which there is no
unbiased estimator must have variance diverging to infinity if the bias
decreases to zero. Given that one cannot eliminate the bias entirely without
its variance exploding, the bias reduction may do the estimator either harm
or good. (See Hirano and Porter (2012) for a recent result for
nondifferentiable parameters.)

Many early researches on estimation of a nonregular parameter considered a
parametric model and focused on finite sample optimality properties. For
example, estimation of a normal mean under bound restrictions or order
restrictions has been studied, among many others, by Lovell and Prescott
(1970), Casella and Strawderman (1981), Bickel (1981), Moors (1981), and
more recently van Eeden and Zidek (2004). Closer to this paper are
researches by Blumenthal and Cohen (1968a,b) who studied estimation of $\max
\{\beta _{1},\beta _{2}\},$ when i.i.d. observations from a location family
of symmetric distributions or normal distributions are available. On the
other hand, the notion of asymptotic efficient estimation through the
convolution theorem and the local asymptotic minimax theorem initiated by Haj%
\'{e}k (1972) and Le Cam (1979) has mostly focused on regular parameters,
and in many cases, resulted in regular estimators as optimal estimators.
Hence the classical theory of semiparametric estimation widely known and
well summarized in monographs such as Bickel, Klassen, Ritov, and Wellner
(1993) and in later sections of van der Vaart and Wellner (1996) (Sections
3.10-3.11, pp. 401-422) does not directly apply to the problem of estimation
of $\theta =(f\circ g)(\mathbf{\beta })$. This paper attempts to fill this
gap from the perspective of local asymptotic minimax estimation.

This paper finds that for the class of nonregular parameters of the form (%
\ref{th}), we can extend the existing theory of local asymptotic minimax
estimation and construct a reasonable class of optimal estimators that are
nonregular in general and asymptotically biased. The class of optimal
estimators take the form of a plug-in estimator with semiparametrically
efficient estimator of $\mathbf{\beta }$ except that it involves an additive
bias-adjustment term which can be computed using simulations.

To deal with nondifferentiability, this paper first focuses on the special
case where $f$ is an identity, and utilizes the approach of generalized
convolution theorem in van der Vaart (1989) to establish the local
asymptotic minimax risk bound for the parameter $\theta $. However, such a
risk bound is hard to use in our set-up where $f$ or $g$ is potentially
asymmetric, because the risk bound involves minimization of the risk over
the distributions of \textquotedblleft noise\textquotedblright\ in the
convolution theorem. This paper proposes a local asymptotic minimax
decision of a simple form: 
\begin{equation*}
g(\mathbf{\hat{\beta}})+\hat{c}/\sqrt{n},
\end{equation*}%
where $\mathbf{\hat{\beta}}$ is a semiparametrically efficient estimator of $%
\mathbf{\beta }$ and $\hat{c}$ is a bias adjustment term that can be
computed through simulations.

Next, extension to the case where $f$ is continuous piecewise linear with a
single kink point is done. Thus, an estimator of the form%
\begin{equation}
\hat{\theta}_{mx}\equiv f\left( g(\mathbf{\hat{\beta}})+\frac{\hat{c}}{\sqrt{%
n}}\right) ,  \label{oe}
\end{equation}%
with appropriate bias adjustment term $\hat{c}$, is shown to be local
asymptotic minimax. In several situations, the bias adjustment term $\hat{c}$
can be set to zero. In particular, when $\theta =\mathbf{s}^{\top }\mathbf{%
\beta }$, for some known vector $\mathbf{s}\in \mathbf{R}^{d}$, so that $%
\theta $ is a regular parameter, the bias adjustment term can be set to be
zero, and an optimal estimator in (\ref{oe}) is reduced to $\mathbf{s}^{\top
}\mathbf{\hat{\beta}}$ which is a semiparametric efficient estimator of $%
\theta =\mathbf{s}^{\top }\mathbf{\beta }$. This confirms the continuity of
this paper's approach with the standard method of semiparametric efficiency.

This paper offers results from a small sample simulation study for the case
of $\theta =\max \{\beta _{1},\beta _{2}\}$. This paper compares the method
with two alternative bias reduction methods:\ fixed bias reduction method
and a selective bias reduction method. The method of local asymptotic
minimax estimation shows relatively robust performance in terms of the
finite sample risk.

The next section defines the scope of the paper by introducing
nondifferentiable transforms that this paper focuses on. The section also
introduces regularity conditions for probabilities that identify $\mathbf{%
\beta }$. Section 3 investigates optimal decisions based on the local
asymptotic maximal risks. Section 4 presents and discusses Monte Carlo
simulation results. All the mathematical proofs are relegated to the
Appendix.

\section{Nondifferentiable Transforms of a Regular Parameter}

In this section, we present the details of the set-up in this paper. We
introduce some notation. Let $\mathbb{N}$ be the collection of natural
numbers. Let $\mathbf{1}_{d}$ be a $d\times 1$ vector of ones with $d\geq 2$%
. For a vector $\mathbf{x}\in \mathbf{R}^{d}$ and a scalar $c$, we simply
write $\mathbf{x}+c=\mathbf{x}+c\mathbf{1}_{d}$, or write $\mathbf{x}=c$
instead of $\mathbf{x}=c\mathbf{1}_{d}$. We define $S_{1}\equiv \{\mathbf{x}%
\in \mathbf{R}^{d}:\mathbf{x}^{\top }\mathbf{1}_{d}=1\}$, where the notation 
$\equiv $ indicates definition. For $\mathbf{x}\in \mathbf{R}^{d}$, the
notation $\max (\mathbf{x})$ (or $\min (\mathbf{x})$) means the maximum (or
the minimum) over the entries of the vector $\mathbf{x}$. When $x_{1},\cdot
\cdot \cdot ,x_{n}$ are scalars, we also use the notations $\max
\{x_{1},\cdot \cdot \cdot ,x_{n}\}$ and $\min \{x_{1},\cdot \cdot \cdot
,x_{n}\}$ whose meanings are obvious. We let $\mathbf{\bar{R}}=[-\infty
,\infty ]$ and view it as a two-point compactification of $\mathbf{R}$, and
let $\mathbf{\bar{R}}^{d}$ be the product of its $d$ copies, so that $%
\mathbf{\bar{R}}^{d}$ itself is a compactification of $\mathbf{R}^{d}$.
(e.g. Dudley (2002), p.74.) We follow the convention to set $\infty \cdot
0=0 $ and $(-\infty )\cdot 0=0$. A supremum and an infimum of a nonnegative
map over an empty set are set to be 0 and $\infty $ respectively.

As for the parameter of interest $\theta $, this paper assumes that%
\begin{equation}
\theta =(f\circ g)(\mathbf{\beta }),  \label{param}
\end{equation}%
where $\mathbf{\beta }\in \mathbf{R}^{d}$ is a regular parameter (the
meaning of regularity for $\mathbf{\beta }$ is clarified in Assumption 2
below), and $g:\mathbf{R}^{d}\rightarrow \mathbf{R}$ and $f:\mathbf{R}%
\rightarrow \mathbf{R}$ satisfy the following assumptions.\bigskip

\noindent \textsc{Assumption 1:} (i) The map $g:\mathbf{R}^{d}\rightarrow 
\mathbf{R}$ is Lipschitz continuous, and satisfies the following.

(a) (Translation Equivariance) For each $c\in \mathbf{R}$ and $\mathbf{x}\in 
\mathbf{R}^{d}$, $g(\mathbf{x}+c)=g(\mathbf{x})+c.$

(b) (Scale Equivariance) For each $u\geq 0$ and $\mathbf{x}\in \mathbf{R}%
^{d},$ $g(u\mathbf{x})=ug(\mathbf{x}).$

(c) (Directional Derivatives) For each $\mathbf{z}\in \mathbf{R}^{d}$ and $%
\mathbf{x}\in \mathbf{R}^{d}$,%
\begin{equation*}
\tilde{g}(\mathbf{x};\mathbf{z})\equiv \lim_{t\downarrow 0}t^{-1}\left(
g\left( \mathbf{x}+t\mathbf{z}\right) -g\left( \mathbf{x}\right) \right)
\end{equation*}%
exists.

\noindent (ii) The map $f:\mathbf{R}\rightarrow \mathbf{R}$ is continuous,
piecewise linear with one kink at a point (i.e., one point of nonlinearity)
in $\mathbf{R}$.\bigskip

We collect here the properties of the directional derivative $\tilde{g}(%
\mathbf{x};\mathbf{z})$ in (c) of the translation-scale equivariant and
Lipschitz continuous map $g$.\bigskip

\noindent \textsc{Lemma 1:} (i) For\textit{\ }each $\mathbf{z}\in \mathbf{R}%
^{d},$\textit{\ }$\mathbf{x}\in \mathbf{R}^{d}$, $c\in \mathbf{R}$, and $%
u\geq 0$, the following properties are satisfied:

(a) $\tilde{g}(\mathbf{0};\mathbf{z})=g(\mathbf{z}).$

(b)$\text{ }\tilde{g}(\mathbf{x}+c;\mathbf{z})=\tilde{g}(\mathbf{x};\mathbf{z%
}).$

(c) $\tilde{g}(\mathbf{x};\mathbf{z}+c)=\tilde{g}(\mathbf{x};\mathbf{z})+c.$

(d) $\tilde{g}(u\mathbf{x};u\mathbf{z})=u\tilde{g}(\mathbf{x};\mathbf{z})=%
\tilde{g}(\mathbf{x};u\mathbf{z}).$

\noindent (ii) For each $\mathbf{x}\in \mathbf{R}^{d}$, $\tilde{g}(\mathbf{x}%
;\mathbf{z})$ is Lipschitz continuous in $\mathbf{z}\in \mathbf{R}^{d}$.

\noindent (iii) For\textit{\ }each $\mathbf{x}\in \mathbf{R}^{d}$, the
convergence in the definition of the directional derivative in Assumption
1(i)(c) is uniform over $\mathbf{z}$ in any bounded subset of $\mathbf{R}%
^{d} $.\bigskip

Assumption 1 essentially defines the scope of this paper. Some examples of $%
g $ are as follows.\bigskip

\noindent \textsc{Examples 3:} (a) $g(\mathbf{x})=\mathbf{s}^{\top }\mathbf{x%
},$ where $\mathbf{s}\in S_{1}$.

\noindent (b) $g(\mathbf{x})=\max (\mathbf{x})$ or $g(\mathbf{x})=\min (%
\mathbf{x})$.

\noindent (c) $g(\mathbf{x})=\max \{\min (\mathbf{x}_{1}),\mathbf{x}_{2}\}$, 
$g(\mathbf{x})=\max (\mathbf{x}_{1})+\max (\mathbf{x}_{2}),$ $g(\mathbf{x}%
)=\min (\mathbf{x}_{1})+\min (\mathbf{x}_{2}),$ $g(\mathbf{x})=\max (\mathbf{%
x}_{1})+\min (\mathbf{x}_{2}),$ or $g(\mathbf{x})=\max (\mathbf{x}_{1})+%
\mathbf{s}^{\top }\mathbf{x}$ with $\mathbf{s}\in S_{1}$, where $\mathbf{x}%
_{1}$ and $\mathbf{x}_{2}$ are subvectors of $\mathbf{x}$. $\blacksquare $%
\bigskip

One might ask whether the representation of parameter $\theta $ as a
composition map $f\circ g$ of $\mathbf{\beta }$ in (\ref{param}) is unique.
The following lemma gives an affirmative answer.\bigskip

\noindent \textsc{Lemma 2:} \textit{Suppose that} $f_{1}$ \textit{and} $%
f_{2} $ \textit{are }$\mathbf{R}$\textit{-valued maps on }$\mathbf{R}$%
\textit{\ that are non-constant on }$\mathbf{R}$, \textit{and }$g_{1}$ 
\textit{and} $g_{2}$ \textit{satisfy} Assumption 1(i). \textit{If} $%
f_{1}\circ g_{1}=f_{2}\circ g_{2},$ \textit{we have}%
\begin{equation*}
f_{1}=f_{2}\text{ \textit{and} }g_{1}=g_{2}\text{\textit{.}}
\end{equation*}%
\bigskip

As we shall see later, the local asymptotic minimax risk bound and the
optimal estimators involve the maps $f$ and $g$. The uniqueness result of
Lemma 2 removes ambiguity that could potentially arise when $\theta $ had
multiple equivalent representations with different maps $f$ and $g$.

We introduce briefly conditions for probabilities that identify $\mathbf{%
\beta }$, in a manner adapted from van der Vaart (1991) and van der Vaart
and Wellner (1996) (see Section 3.11, pp. 412-422.) Let $\mathcal{P}\equiv
\{P_{\alpha }:\alpha \in \mathcal{A}\}$ be a family of distributions on a
measurable space $(\mathcal{X},\mathcal{G})$ indexed by $\alpha \in \mathcal{%
A}$, where the set $\mathcal{A}$ is a nonempty open subset of a Euclidean
space or more generally a complete metric space.

We assume that we have i.i.d. draws $Y_{1},\cdot \cdot \cdot ,Y_{n}$ from $%
P_{\alpha _{0}}\in \mathcal{P}$ for some $\alpha _{0}\in \mathcal{A}$, so
that $\mathbf{X}_{n}\equiv (Y_{1},\cdot \cdot \cdot ,Y_{n})$ is distributed
as $P_{\alpha _{0}}^{n}$. Let $\mathcal{P}(P_{\alpha _{0}})$ be the
collection of maps $t\rightarrow P_{\alpha _{t}}\ $such that for some $h\in
L_{2}(P_{\alpha _{0}})$,%
\begin{equation}
\int \left\{ \frac{1}{t}\left( dP_{\alpha _{t}}^{1/2}-dP_{\alpha
_{0}}^{1/2}\right) -\frac{1}{2}hdP_{\alpha _{0}}^{1/2}\right\}
^{2}\rightarrow 0,\text{ as\ }t\rightarrow 0.  \label{cv}
\end{equation}%
When this convergence holds, we say that $P_{\alpha _{t}}$ is \textit{%
differentiable in quadratic mean} to $P_{\alpha _{0}}$, call $h\in
L_{2}(P_{\alpha _{0}})$ a \textit{score function} associated with this
convergence, and call the set of all such $h$'s a \textit{tangent set},
denoting it by $T(P_{\alpha _{0}}).$ We assume that the tangent set is a
linear subspace of $L_{2}(P_{\alpha _{0}})$. Taking $\langle \cdot ,\cdot
\rangle $ to be the usual inner product in $L_{2}(P_{\alpha _{0}})$, we
write $H\equiv T(P_{\alpha _{0}})$ and view $(H,\langle \cdot ,\cdot \rangle
)$ as a subspace of a separable Hilbert space, with $\bar{H}$ denoting its
completion. For each $h\in H,$ $n\in \mathbb{N}$, and $\lambda _{h}\in 
\mathcal{A},$ let $P_{\alpha _{0}+\lambda _{h}/\sqrt{n}}$ be probabilities
converging to $P_{\alpha _{0}}$ (as in (\ref{cv})) as $n\rightarrow \infty $
having $h$ as its associated score. We simply write $P_{n,h}=P_{\alpha
_{0}+\lambda _{h}/\sqrt{n}}^{n}$ and consider sequences of such
probabilities $\{P_{n,h}\}_{n\geq 1}$ indexed by $h\in H$. (See van der
Vaart (1991) and van der Vaart and Wellner (1996), Section 3.11 for
details.) The collection $\mathcal{E}_{n}\equiv (\mathcal{X}_{n},\mathcal{G}%
_{n},P_{n,h};h\in H)$ constitutes a sequence of statistical experiments for $%
\mathbf{\beta }$.

Due to differentiability in quadratic mean and i.i.d. assumption, the
collection $\mathcal{E}_{n}$ satisfies \textit{local asymptotic normality}
(LAN), that is, for any $h\in H,$%
\begin{equation*}
\log \frac{dP_{n,h}}{dP_{n,0}}=\zeta _{n}(h)-\frac{1}{2}\langle h,h\rangle ,
\end{equation*}%
where for any $h,h^{\prime }\in H$, $[\zeta _{n}(h),\zeta _{n}(h^{\prime })]%
\overset{d}{\rightarrow }[\zeta (h),\zeta (h^{\prime })]$, under $%
\{P_{n,0}\} $\ and $\zeta (\cdot )$ is a centered Gaussian process on $H\ $%
with covariance function $\mathbf{E}[\zeta (h_{1})\zeta (h_{2})]=\langle
h_{1},h_{2}\rangle $. Note that we require here the joint convergence of $%
\zeta _{n}(h)$ and $\zeta _{n}(h^{\prime })$ for each pair $(h,h^{\prime })$%
. This joint convergence is used to derive a modified version of LAN (Lemma
A4 in the appendix) which is used to derive the local asymptotic minimax
risk. The joint convergence can be seen to hold e.g. from the proof of Lemma
3.10.11 of van der Vaart and Wellner (1996), p.406.

The LAN property reduces the decision problem to one in which an optimal
decision is sought under a single Gaussian shift experiment $\mathcal{E}=(%
\mathcal{X},\mathcal{G},P_{h};h\in H),$ where $P_{h}$ is such that $\log
dP_{h}/dP_{0}=\zeta (h)-\frac{1}{2}\langle h,h\rangle .$

The parameter $\mathbf{\beta }$ is represented as a functional $\mathbf{%
\beta }:\mathcal{P}\rightarrow \mathbf{R}^{d}$. From here on, we simply
write for each $\alpha \in \mathcal{A},$ $\mathbf{\beta }_{n}(h)=\mathbf{%
\beta }(P_{\alpha _{0}+\lambda _{h}/\sqrt{n}})$ and regard $\mathbf{\beta }%
_{n}(\cdot )$ as an $\mathbf{R}^{d}$-valued map on $H$.\bigskip

\noindent \textsc{Assumption 2:} (Regular Parameter) There exists a
continuous linear $\mathbf{R}^{d}$-valued map, $\mathbf{\dot{\beta}}$, on $H$
such that for any $h\in H,$%
\begin{equation*}
\sqrt{n}(\mathbf{\beta }_{n}(h)-\mathbf{\beta }_{n}(0))\rightarrow \mathbf{%
\dot{\beta}}(h),
\end{equation*}%
as $n\rightarrow \infty .$\bigskip

Assumption 2 requires that $\mathbf{\beta }\ $be \textit{regular} in the
sense of van der Vaart \ and Wellner (1996, Section 3.11). The map $\mathbf{%
\dot{\beta}}$ in Assumption 2 is associated with the semiparametric
efficiency bound of $\mathbf{\beta }$. For each $\mathbf{b}\in \mathbf{R}%
^{d} $, $\mathbf{b}^{\top }\mathbf{\dot{\beta}}(\cdot )$ defines a
continuous linear functional on $H$, and hence there exists $\dot{\beta}_{%
\mathbf{b}}^{\ast }\in \bar{H}$ such that $\mathbf{b}^{\top }\mathbf{\dot{%
\beta}}(h)=\langle \dot{\beta}_{\mathbf{b}}^{\ast },h\rangle ,$ $h\in H$.
Then for any $\mathbf{b}\in \mathbf{R}^{d}$, $||\dot{\beta}_{\mathbf{b}%
}^{\ast }||^{2} $ represents the asymptotic variance bound of the parameter $%
\mathbf{b}^{\top }\mathbf{\beta }$. The map $\dot{\beta}_{\mathbf{b}}^{\ast
} $ is called an \textit{efficient influence function} for $\mathbf{b}^{\top
}\mathbf{\beta }$ in the literature (e.g. van der Vaart (1991)). Let $%
\mathbf{e}_{m}$ be a $d\times 1$ vector whose $m$-th entry is one and the
other entries are zero, and let $\Sigma $ be a $d\times d$ matrix whose $%
(m,k)$-th entry is given by $\langle \dot{\beta}_{\mathbf{e}_{m}}^{\ast },%
\dot{\beta}_{\mathbf{e}_{k}}^{\ast }\rangle $. As for $\Sigma $, we assume
the following:\bigskip

\noindent \textsc{Assumption 3}: $\Sigma $ is invertible.\bigskip

\noindent The inverse of matrix $\Sigma $ is called the semiparametric
efficiency bound for $\mathbf{\beta }.$ In particular, Assumption 3 requires
that there is no redundancy among the entries of $\mathbf{\beta }$, i.e.,
one entry of $\mathbf{\beta }$ is not defined as a linear combination of the
other entries.

\section{Local Asymptotic Minimax Estimators}

\subsection{Loss Functions}

For a decision $d\in \mathbf{R}$ and the object of interest $\theta \in 
\mathbf{R}$, we consider the following form of a loss function:%
\begin{equation}
L\left( d,\theta \right) =\tau (|d-\theta |),  \label{LB}
\end{equation}%
where $\tau :\mathbf{R}\rightarrow \mathbf{R}$ is a map that satisfies the
following assumption.\bigskip

\noindent \textsc{Assumption 4}: (i) $\tau (\cdot )$ is increasing and convex on $%
[0,\infty )$, $\tau (0)=0$, and there exists $\bar{\tau}\in (0,\infty ]$ such%
$\ $that $\tau ^{-1}([0,y])$ is bounded in $[0,\infty )$ for all $0<y<\bar{%
\tau}$.

\noindent (ii) For each $M>0$, there exists $C_{M}>0$ such that for all $%
x,y\in \mathbf{R}$, 
\begin{equation}
|\tau _{M}(x)-\tau _{M}(y)|\ \leq \ C_{M}|x-y|,  \label{Lip}
\end{equation}%
where $\tau _{M}(\cdot )=\min \{\tau (\cdot ),M\}$.\bigskip

The smoothness condition in (\ref{Lip}) is weaker than requiring $\tau $ to be Lipschitz
continuous. For example, the squared loss function $\tau (x)=x^{2}$
satisfies this condition, but not Lipschitz continuity. While Assumption 4
is satisfied by many loss functions, it excludes the hypothesis testing type
loss function $\tau (|d-\theta |)=1\{|d-\theta |>c\}$, $c\in \mathbf{R}$.
From here on, we identify $\tau $ and $\tau _{M}$ as their continuous
extensions to $(-\infty ,\infty ]$.

The following lemma establishes a lower bound for the local asymptotic
minimax risk when $f$ is an identity. Let for each $b\in \lbrack 0,\infty )$
and $n\geq 1,$%
\begin{equation*}
H_{n,b}\equiv \left\{ h\in H:||\mathbf{\beta }_{n}(h)-\mathbf{\beta }%
_{n}(0)||\leq b/\sqrt{n}\right\} .
\end{equation*}%
The set $H_{n,b}$ collects those $h$'s in $H$ at which $\mathbf{\beta }%
_{n}(h)$ lies locally around $\mathbf{\beta }_{n}(0)$. (Confining our
attention to $h\in H_{n,b}$ enables us to control the convergence in
Assumption 2 uniformly over $h$ in $H_{n,b}$.)\bigskip

\noindent \textsc{Lemma 3}: \textit{Suppose that Assumptions 1-4 hold and
that }$f$ \textit{is an identity. Then for any sequence of estimators }$\hat{%
\theta}$,%
\begin{eqnarray*}
&&\sup_{b\in \lbrack 0,\infty )}\ \underset{n\rightarrow \infty }{\text{%
liminf}}\sup_{h\in H_{n,b}}\mathbf{E}_{h}\left[ \tau (|\sqrt{n}\{\hat{\theta}%
-g(\mathbf{\beta }_{n}(h)\}|)\right] \\
&\geq &\inf_{F\in \mathcal{F}}\sup_{\mathbf{r}\in \mathbf{R}^{d}}\int 
\mathbf{E}\left[ \tau (|\tilde{g}_{0}(Z+\mathbf{r})-\tilde{g}_{0}(\mathbf{r}%
)+w)|)\right] dF(w),
\end{eqnarray*}%
\textit{where} $\mathbf{\beta }_{0}\equiv \mathbf{\beta }(P_{\alpha _{0}})$, 
$\tilde{g}_{0}(\mathbf{r})\equiv \tilde{g}(\mathbf{\beta }_{0};\mathbf{r})$, 
$\mathbf{E}_{h}$ \textit{denotes the expectation under }$P_{n,h}$\textit{,
and }$\mathcal{F}$\textit{\ denotes the collection of probability measures
on the Borel }$\sigma $\textit{-field of }$\mathbf{R}$\textit{.}\bigskip

The lower bound in Lemma 3 involves the directional derivatives $\tilde{g}%
(\cdot ;\cdot )$ of $g$. Typically computation of directional derivatives is
straightforward in many examples. (However, the practical procedure of
optimal estimation proposed in this paper does not require an explicit
computation of the directional derivatives, as we shall see after Assumption
5.)\bigskip

\noindent \textsc{Examples 4:} (a) Suppose that $g(\mathbf{x})=\mathbf{s}%
^{\top }\mathbf{x}$, $\mathbf{s}\in S_{1}$. Then obviously, $\tilde{g}_{0}(%
\mathbf{z})=\mathbf{s}^{\top }\mathbf{z},$ and the risk lower bound in Lemma
3 becomes%
\begin{equation*}
\inf_{F\in \mathcal{F}}\int \mathbf{E}\left[ \tau (|\mathbf{s}^{\top }Z+w|)%
\right] dF(w)\geq \mathbf{E}\left[ \tau (|\mathbf{s}^{\top }Z|)\right] ,
\end{equation*}%
the last inequality following from Anderson's Lemma.

\noindent (b) Suppose that $g(\mathbf{x})=\max \{x_{1},x_{2}\}$. Then%
\begin{equation*}
\tilde{g}_{0}(\mathbf{z})=\left\{ 
\begin{array}{c}
z_{1}\text{, if }\beta _{0,1}>\beta _{0,2} \\ 
z_{2},\text{ if }\beta _{0,1}<\beta _{0,2} \\ 
\max \{z_{1},z_{2}\}\text{, if }\beta _{0,1}=\beta _{0,2},%
\end{array}%
\right.
\end{equation*}%
where $\beta _{0,1}$ and $\beta _{0,2},$ and $z_{1}$ and $z_{2}$ are the
first and the second entries of $\mathbf{\beta }_{0}$ and $\mathbf{z}$
respectively. $\blacksquare $\bigskip

The lower bound in Lemma 3 is obtained by using a version of a generalized
convolution theorem in van der Vaart (1989) which is adapted to the current
set-up. The main difficulty with using Lemma 3 is that the supremum over $%
\mathbf{r}\in \mathcal{S}$ and the infimum over $F\in \mathcal{F}$ do not
have an explicit solution in general. Hence this paper considers simulating
the lower bound in Lemma 3 by using random draws from a distribution
approximating that of $Z$. The main obstacle in this approach is that the
risk lower bound involves infimum over an infinite dimensional space $%
\mathcal{F}$.

We now simplify the risk lower bound. By Jensen's inequality,
\begin{eqnarray*}
	&& \inf_{F\in \mathcal{F}}\sup_{\mathbf{r}\in \mathbf{R}^{d}}\int 
	\mathbf{E}\left[ \tau (|\tilde{g}_{0}(Z+\mathbf{r})-\tilde{g}_{0}(\mathbf{r}%
	)+w|)\right] dF(w)\\
	&\ge& \inf_{F\in \mathcal{F}}\sup_{\mathbf{r}\in \mathbf{R}^{d}}
	\mathbf{E}\left[ \tau \left(\left|\tilde{g}_{0}(Z+\mathbf{r})-\tilde{g}_{0}(\mathbf{r}%
	)+\int w dF(w)\right|\right)\right].
\end{eqnarray*}
Thus we obtain the following theorem.\bigskip

\noindent \textsc{Theorem 1}: \textit{Suppose that Assumptions 1-4 hold and
that }$f$\textit{\ is an identity. Then for any sequence of estimators }$%
\hat{\theta}$,%
\begin{equation*}
\sup_{b\in \lbrack 0,\infty )}\ \underset{n\rightarrow \infty }{\text{liminf}%
}\sup_{h\in H_{n,b}}\mathbf{E}_{h}\left[ \tau (|\sqrt{n}\{\hat{\theta}-g(%
\mathbf{\beta }_{n}(h))\}|)\right] \geq \inf_{c\in \mathbf{R}}B(c;1),
\end{equation*}%
\textit{where for }$c\in \mathbf{R}$, \textit{and} \textit{any }$a\geq 0,$%
\begin{equation*}
B(c;a)\equiv \sup_{\mathbf{r}\in \mathbf{R}^{d}}\mathbf{E}\left[ \tau (a|%
\tilde{g}_{0}(Z+\mathbf{r})-\tilde{g}_{0}(\mathbf{r})+c)|)\right] .
\end{equation*}%
\bigskip

The main feature of the lower bound in Theorem 1 is that it involves infimum
over a \textit{single-dimensional} space $\mathbf{R}$ in its risk bound.
This simpler form now makes it feasible to simulate the lower bound for the
risk.

This paper proposes a method of constructing a local asymptotic minimax
estimator as follows. Suppose that we are given a consistent estimator $\hat{%
\Sigma}$ of $\Sigma $ and a semiparametrically efficient estimator $\mathbf{%
\hat{\beta}}$ of $\mathbf{\beta }$ which satisfy the following assumptions.
(See Bickel, Klaasen, Ritov, and Wellner (1993) for semiparametric efficient
estimators from various models.)\bigskip

\noindent \textsc{Assumption 5:} (i) For each $\varepsilon >0$, there exists 
$M>0$ such that 
\begin{equation*}
\underset{n\rightarrow \infty }{\text{limsup}}\sup_{h\in H}P_{n,h}\{\sqrt{n}%
||\hat{\Sigma}-\Sigma ||>M\}<\varepsilon \text{.}
\end{equation*}

\noindent (ii) For each $t\in \mathbf{R}^{d}$,%
\begin{equation*}
\underset{n\rightarrow \infty }{\text{limsup}}\sup_{h\in H}\left\vert
P_{n,h}\{\sqrt{n}(\mathbf{\hat{\beta}}-\mathbf{\beta }_{n}(h))\leq
t\}-P\{Z\leq t\}\right\vert =0,
\end{equation*}%
as $n\rightarrow \infty $.\bigskip

Assumption 5 imposes $\sqrt{n}$-consistency of $\hat{\Sigma}$ and
convergence in distribution of $\sqrt{n}(\mathbf{\hat{\beta}}-\mathbf{\beta }%
_{n}(h)),$ both uniform over $h\in H$. The uniform convergence can be proved
through the central limit theorem uniform in $h\in H$. Under regularity
conditions, the uniform central limit theorem of a sum of i.i.d. random
variables follows from a Berry-Esseen bound, as long as the third moment of
the random variable is bounded uniformly in $h\in H.$

For a fixed large $M_{1}>0,$ we define%
\begin{equation}
\hat{\theta}_{mx}\equiv g(\mathbf{\hat{\beta})}+\frac{\hat{c}_{M_{1}}}{\sqrt{%
n}},  \label{mx}
\end{equation}%
where $\hat{c}_{M_{1}}$ is a bias adjustment term constructed from the
simulations of the risk lower bound in Theorem 1, as we explain now. (Note
that $\hat{\theta}_{mx}$ depends on $M_{1}$ in general though the dependence
is suppressed from notation.)

To simulate the risk lower bound in Theorem 1, we first draw $\{\mathbf{\xi }%
_{i}\}_{i=1}^{L}$ i.i.d. from $N(0,I_{d})$. Since $\tilde{g}_{0}(\cdot )$
depends on $\mathbf{\beta }_{0}$ that is unknown to the researcher, we first
construct a consistent estimator of $\tilde{g}_{0}(\cdot ).$ Take a sequence 
$\varepsilon _{n}\rightarrow 0$ such that $\sqrt{n}\varepsilon
_{n}\rightarrow \infty $ as $n\rightarrow \infty $. Examples of $\varepsilon
_{n}$ are $\varepsilon _{n}=n^{-1/3}$ or $\varepsilon _{n}=n^{-1/2}\log n$.
Observe that $\tilde{g}_{0}(\mathbf{z}),$ $\mathbf{z\in R}^{d}$, is
approximated by 
\begin{equation*}
\varepsilon _{n}^{-1}\left( g\left( \varepsilon _{n}\mathbf{z+\beta }%
_{0}\right) -g(\mathbf{\beta }_{0}\mathbf{)}\right) =g\left( \mathbf{z+}%
\varepsilon _{n}^{-1}(\mathbf{\beta }_{0}-g(\mathbf{\beta }_{0}))\right) ,
\end{equation*}%
as $n\rightarrow \infty $. Hence we define%
\begin{equation*}
\hat{g}_{n}(\mathbf{z})\equiv g\left( \mathbf{z+}\varepsilon _{n}^{-1}(%
\mathbf{\hat{\beta}}-g(\mathbf{\hat{\beta}}))\right) .
\end{equation*}%
Then it is not hard to see that $\hat{g}_{n}(\mathbf{z})$ is consistent for $%
\tilde{g}_{0}(\mathbf{z})$. Thus, we consider the following: for any $a\geq
0,$%
\begin{equation*}
\hat{B}_{M_{1}}(c;a)\equiv \sup_{\mathbf{r}\in \lbrack -M_{1},M_{1}]^{d}}%
\frac{1}{L}\sum_{i=1}^{L}\tau _{M_{1}}\left( a\left\vert \hat{g}_{n}(\hat{%
\Sigma}^{1/2}\mathbf{\xi }_{i}+\mathbf{r})-\hat{g}_{n}(\mathbf{r}%
)+c\right\vert \right) \text{.}
\end{equation*}%
Then we define%
\begin{equation}
\hat{c}_{M_{1}}(a)\equiv \sup \hat{E}_{M_{1}}(a) ,  \label{cmx}
\end{equation}%
where, with $\eta _{n,L}\rightarrow 0$ as $n,L\rightarrow \infty $, $\eta
_{n,L}\varepsilon _{n}\sqrt{n}\rightarrow \infty $ as $n\rightarrow \infty $
and $\eta _{n,L}\sqrt{L}\rightarrow \infty $ as $L\rightarrow \infty $,%
\begin{equation*}
\hat{E}_{M_{1}}(a)\equiv \left\{ c\in \lbrack -M_{1},M_{1}]:\hat{B}%
_{M_{1}}(c;a)\leq \inf_{c_{1}\in \lbrack -M_{1},M_{1}]}\hat{B}%
_{M_{1}}(c_{1};a)+\eta _{n,L}\right\} .
\end{equation*}%
The formulation of $\hat{c}_{M_{1}}(a)$ in (\ref{cmx}) is designed to yield
an unambiguous determination of a minimizer of $\hat{B}_{M_{1}}(c;a)$ (up to
a small number $\eta _{n,L}$) over $c\in \lbrack -M_{1},M_{1}]$, even when
the minimizer of its population version $B(c;a)$ over $c\in \lbrack
-M_{1},M_{1}]$ turns out to be non-unique.

Now, as for the bias adjustment term $\hat{c}_{M_{1}}$ in (\ref{mx}), we
take $\hat{c}_{M_{1}}=\hat{c}_{M_{1}}(1)$. The following theorem affirms
that $\hat{\theta}_{mx}$ is local asymptotic minimax for $\theta =g(\mathbf{%
\beta })$. (For technical facility, we follow a suggestion by Strasser
(1985) (p.440) and consider a truncated loss: $\tau _{M}(\cdot )=\min \{\tau
(\cdot ),M\}$ for large $M.$)\bigskip

\noindent \textsc{Theorem 2:} \textit{Suppose that the conditions of Theorem
1 and Assumption 5 hold. Then for any }$M>0$\textit{\ and any }$M_{1}\geq M$%
\textit{\ that constitutes }$\hat{c}_{M_{1}}$,\textit{\ }%
\begin{equation*}
\sup_{b\in \lbrack 0,\infty )}\ \underset{n\rightarrow \infty }{\text{limsup}%
}\sup_{h\in H_{n,b}}\mathbf{E}_{h}\left[ \tau _{M}(|\sqrt{n}\{\hat{\theta}%
_{mx}-g(\mathbf{\beta }_{n}(h))\}|)\right] \leq \inf_{c\in \mathbf{R}}B(c;1).
\end{equation*}%
\bigskip

Recall that the candidate estimators considered in Theorem 1 were not
restricted to plug-in estimators with an additive bias adjustment term. As
standard in the literature of local asymptotic minimax estimation, the
candidate estimators are any sequences of measurable functions of
observations including both regular and nonregular estimators. The main
thrust of Theorem 2 is the finding that it is sufficient for local
asymptotic minimax estimation to consider a plug-in estimator using a
semiparametrically efficient estimator of $\mathbf{\beta }$ with an additive
bias adjustment term as in (\ref{mx}). It remains to find optimal bias
adjustment, which can be done using the simulation method proposed earlier.

We now extend the result to the case where $f$ is not an identity map, but a
continuous piecewise linear map with a single kink point $\bar{x}\in \mathbf{%
R}$. For concreteness, suppose that for all $x\in \mathbf{R}$,%
\begin{equation*}
f(x)=\left\{ 
\begin{array}{c}
a_{1}(x-\bar{x})+f(\bar{x})\text{, if }x\geq \bar{x} \\ 
a_{2}(x-\bar{x})+f(\bar{x})\text{, if }x<\bar{x}%
\end{array}%
\right.
\end{equation*}%
for $a_{1},a_{2}\in \mathbf{R}$. Let%
\begin{equation*}
s\equiv \left\{ 
\begin{array}{c}
|a_{1}|\text{,} \\ 
|a_{2}|\text{,} \\ 
\max \left\{ |a_{1}|,|a_{2}|\right\} \text{,}%
\end{array}%
\begin{array}{c}
\text{if }g\left( \mathbf{\beta }_{0}\right) >\bar{x} \\ 
\text{if }g\left( \mathbf{\beta }_{0}\right) <\bar{x} \\ 
\text{if }g\left( \mathbf{\beta }_{0}\right) =\bar{x}%
\end{array}%
\right\} .
\end{equation*}%
Then the following theorem establishes the risk lower bound for the case
where $f$ is not an identity map.\bigskip

\noindent \textsc{Theorem 3:} \textit{Suppose that Assumptions 1-4 hold.
Then for any sequence of estimators }$\hat{\theta}$,%
\begin{equation*}
\sup_{b\in \lbrack 0,\infty )}\ \underset{n\rightarrow \infty }{\text{liminf}%
}\sup_{h\in H_{n,b}}\mathbf{E}_{h}\left[ \tau (|\sqrt{n}\{\hat{\theta}%
-(f\circ g)(\mathbf{\beta }_{n}(h))\}|)\right] \geq \inf_{c\in \mathbf{R}%
}B(c;s).
\end{equation*}%
\bigskip

The bounds in Theorems 1 and 3 involve a bias adjustment term $c^{\ast }$
that minimizes $B(c;s)$ over $c\in \mathbf{R}$. A similar bias adjustment
term appears in Takagi (1994)'s local asymptotic minimax estimation result.
While the bias adjustment term arises here due to asymmetric
nondifferentiable map $f\circ g$ of a regular parameter, it arises in his
paper due to an asymmetric loss function, and the decision problem in this
paper cannot be reduced to his set-up, even if we assume a parametric family
of distributions indexed by an open interval as he does in his paper.

Now let us search for a class of local asymptotic minimax estimators that
achieve the lower bound in Theorem 3. Let 
\begin{equation*}
\hat{s}=\left\{ 
\begin{array}{c}
|a_{1}|\text{,} \\ 
|a_{2}|\text{,} \\ 
\max \left\{ |a_{1}|,|a_{2}|\right\} \text{,}%
\end{array}%
\begin{array}{l}
\text{if }g(\mathbf{\hat{\beta})}>\bar{x}+\varepsilon _{n} \\ 
\text{if }g(\mathbf{\hat{\beta})}<\bar{x}-\varepsilon _{n} \\ 
\text{if }\bar{x}-\varepsilon _{n}\leq g(\mathbf{\hat{\beta})}\leq \bar{x}%
+\varepsilon _{n}%
\end{array}%
\right\} ,
\end{equation*}%
where $\varepsilon _{n}\rightarrow 0$ such that $\sqrt{n}\varepsilon
_{n}\rightarrow \infty $ as $n\rightarrow \infty $. It turns out that an
estimator of the form:%
\begin{equation}
\tilde{\theta}_{mx}\equiv f\left( g(\mathbf{\hat{\beta})}+\frac{\hat{c}%
_{M_{1}}(\hat{s})}{\sqrt{n}}\right) ,  \label{mx2}
\end{equation}%
where $\hat{c}_{M_{1}}(\hat{s})$ is the bias-adjustment term defined in (\ref%
{cmx}) only with $a$ there replaced by $\hat{s}$, is local asymptotic
minimax.\bigskip

\noindent \textsc{Theorem 4:} \textit{Suppose that the conditions of Theorem
3 and Assumption 5 hold. Then},\textit{\ for any }$M>0$\textit{\ and any }$%
M_{1}\geq M$\textit{,}%
\begin{equation*}
\sup_{b\in \lbrack 0,\infty )}\ \underset{n\rightarrow \infty }{\text{limsup}%
}\sup_{h\in H_{n,b}}\mathbf{E}_{h}\left[ \tau _{M}(|\sqrt{n}\{\tilde{\theta}%
_{mx}-(f\circ g)(\mathbf{\beta }_{n}(h))\}|)\right] \leq \inf_{c\in \mathbf{R%
}}B(c;s).
\end{equation*}%
\bigskip

The estimator $\tilde{\theta}_{mx}$ is in general a nonregular estimator
that is asymptotically biased. When $\tau (x)=x^{k}$, $k\geq 1$, we have 
\begin{equation*}
\inf_{c\in \mathbf{R}}B(c;s)=s^{k}\inf_{c\in \mathbf{R}}B(c;1).
\end{equation*}%
Hence it suffices to use $\hat{c}_{M_{1}}(1)$ instead of $\hat{c}_{M_{1}}(%
\hat{s})$ with large $M_{1}$ in this case.

When $g(\mathbf{\beta })=\mathbf{s}^{\top }\mathbf{\beta }$ with $\mathbf{s}%
\in S_{1}$, the risk bound in Theorem 4 becomes%
\begin{equation*}
\inf_{c\in \mathbf{R}}\mathbf{E}\left[ \tau \left( s|\tilde{g}%
_{0}(Z)+c|\right) \right] =\mathbf{E}\left[ \tau \left( s|\mathbf{s}^{\top
}Z|\right) \right] ,
\end{equation*}%
where the equality follows by Anderson's Lemma. In this case, it suffices to
set $\hat{c}_{M_{1}}=0$, because the infimum over $c\in \mathbf{R}$ is
achieved at $c=0$. The minimax decision thus becomes simply%
\begin{equation}
\tilde{\theta}_{mx}=f(\mathbf{\hat{\beta}}^{\top }\mathbf{s}).  \label{sol}
\end{equation}%
This has the following consequences.\bigskip

\noindent \textsc{Examples 5:} (a) When $\theta =\mathbf{\beta }^{\top }%
\mathbf{s}$ for a known vector $\mathbf{s}\in S_{1}$, $\tilde{\theta}_{mx}=%
\mathbf{\hat{\beta}}^{\top }\mathbf{s}$. Therefore, the decision in (\ref%
{sol}) reduces to a semiparametric efficient estimator of $\mathbf{\beta }%
^{\top }\mathbf{s}$.

\noindent (b) When $\theta =\max \{a\mathbf{\beta }^{\top }\mathbf{s}+b,0\}$
for a known vector $\mathbf{s}\in S_{1}$ and known constants $a,b\in \mathbf{%
R}$, $\tilde{\theta}_{mx}=\max \{a\mathbf{\hat{\beta}}^{\top }\mathbf{s}%
+b,0\}.$

\noindent (c) When $\theta =|\beta |$ for a scalar parameter $\beta $, $%
\tilde{\theta}_{mx}=|\hat{\beta}|.$ $\blacksquare $\bigskip

The examples of (b)-(c)\ involve nondifferentiable transform $f$, and hence $%
\tilde{\theta}_{mx}$ as an estimator of $\theta $ is asymptotically biased
in these examples. Nevertheless, the plug-in estimator $\tilde{\theta}_{%
\text{$mx$}}$ that does not require any bias adjustment is local asymptotic
minimax. We provide another example that has the optimal bias adjustment
term equal to zero. This example is motivated by Blumenthal and Cohen
(1968a).\bigskip

\noindent \textsc{Examples 6:} Suppose that $\theta =\max \{\beta _{1},\beta
_{2}\}$, where $\mathbf{\beta }=(\beta _{1},\beta _{2})\in \mathbf{R}^{2}$
is a regular parameter, and the $2\times 2$ matrix $\Sigma $ has identical
diagonal entries equal to $\sigma ^{2}$. (That is, $\beta _{1}$ and $\beta
_{2}$ have the same semiparametric efficiency bound.) We take $\tau
(x)=x^{2} $, i.e., the squared error loss. Then from Example 4(b), the risk
lower bound becomes $\sigma ^{2},$ if $\beta _{0,1}>\beta _{0,2}$ or $\beta
_{0,1}<\beta _{0,2},$ and becomes%
\begin{equation*}
\inf_{c\in \mathbf{R}}\sup_{r\geq 0}\mathbf{E}\left( \max
\{Z_{1}-r,Z_{2}\}-c\right) ^{2},
\end{equation*}%
if $\beta _{0,1}=\beta _{0,2}$, where $Z_{1}$ and $Z_{2}$ denote the first
and second entries of $Z$ respectively.

For each $c\in \mathbf{R}$, $\mathbf{E}\left( \max
\{Z_{1}-r,Z_{2}\}-c\right) ^{2}$ is quasiconvex in $r\geq 0$ so that the
supremum over $r\geq 0$ is achieved at $r=0$ or $r\rightarrow \infty .$ When 
$r=0$, the bound becomes $Var(\max \{Z_{1},Z_{2}\})$ and when $r\rightarrow
\infty $, the bound becomes $Var(Z_{2})$. By (5.10) of Moriguti (1951), we
have $Var(\max \{Z_{1},Z_{2}\})\leq Var(Z_{2})$, so that the local
asymptotic risk bound becomes $Var(Z_{2})=\sigma ^{2}$ with $r=\infty $ and $%
c=0$. Therefore, regardless of $\beta _{0,1}>\beta _{0,2}$, $\beta
_{0,1}<\beta _{0,2}$, or $\beta _{0,1}=\beta _{0,2}$, the risk lower bound
becomes $\sigma ^{2}$ in this case. On the other hand, it is not hard to see
from (A.3) of Blumenthal and Cohen (1968b) that $\tilde{\theta}_{mx}=\max \{%
\hat{\beta}_{1},\hat{\beta}_{2}\}$ (without the bias adjustment term) is
local asymptotic minimax. This result parallels the finding by Blumenthal
and Cohen (1968a) that for squared error loss and observations of two
independent random variables $X_{1}$ and $X_{2}$ from a location family of
symmetric distributions, $\max \{X_{1},X_{2}\}$ is a minimax decision. $%
\blacksquare $

\section{Monte Carlo Simulations}

\subsection{Simulation Designs}

In the simulation study, this paper compares the finite sample risk
performances of the local asymptotic minimax estimator proposed in this
paper with estimators that perform bias reductions in two methods:\ fixed
bias reduction and selective bias reduction.

In this study, we considered the following data generating process. Let $%
\{X_{i}\}_{i=1}^{n}$ be i.i.d. random vectors in $\mathbf{R}^{2}$ where $%
X_{1}\sim N\left( \mathbf{\beta },\Sigma \right) ,$%
\begin{equation}
\mathbf{\beta }=\left[ 
\begin{array}{c}
\beta _{1} \\ 
\beta _{2}%
\end{array}%
\right] =\left[ 
\begin{array}{c}
0 \\ 
\delta _{0}/\sqrt{n}%
\end{array}%
\right] \text{ and }\Sigma =\left[ 
\begin{array}{c}
2 \\ 
1/2%
\end{array}%
\begin{array}{c}
1/2 \\ 
4%
\end{array}%
\right] ,  \label{var2}
\end{equation}%
and $\delta _{0}$ is chosen from grid points in $[-10,10]$. The parameters
of interest are as follows:%
\begin{equation*}
\theta _{1}\equiv f_{1}(g_{1}(\mathbf{\beta }))\text{ and }\theta _{2}\equiv
f_{2}(g_{2}(\mathbf{\beta })),
\end{equation*}%
where 
\begin{eqnarray*}
f_{1}(x) &=&x\text{ and }g_{1}(\mathbf{\beta })=\max \{\beta _{1},\beta
_{2}\},\text{ and} \\
f_{2}(x) &=&\max \{x,0\}\text{ and }g_{2}(\mathbf{\beta })=\beta _{1}.
\end{eqnarray*}%
When $\delta _{0}$ is close to zero, parameters $\theta _{1}$ and $\theta
_{2}$ have $\mathbf{\beta }$ close to the kink point of the
nondifferentiable map. However, when $\delta _{0}$ is away from zero, the
parameters become more like a regular parameter themselves. We take $\mathbf{%
\hat{\beta}}=\frac{1}{n}\sum_{i=1}^{n}X_{i}$ as the estimator of $\mathbf{%
\beta }$. As for the finite sample risk, we adopt the mean squared error:%
\begin{equation*}
\mathbf{E}\left[ (\hat{\theta}_{j}-\theta _{j})^{2}\right] ,\ j=1,2,
\end{equation*}%
where $\hat{\theta}_{j}$ is a candidate estimator for $\theta _{j}$. In the
simulation study, we investigate the finite sample risk profile of decisions
by varying $\delta _{0}$.

We evaluated the risk using Monte Carlo simulations.\ The sample size was
300. The Monte Carlo simulation number was set to be 20,000. The sequence $%
\varepsilon _{n}$ was taken to be $n^{-1/3}$.

We report only the results for the case of $\theta _{1}=f_{1}(g_{1}(\mathbf{%
\beta }))$. The results for the case of $\theta _{2}=f_{2}(g_{2}(\mathbf{%
\beta }))$ were similar and hence omitted.
\begin{figure}[tbph]
\caption{Comparison of the Local Asymptotic Minimax Estimators with
Estimators Obtained through Other Bias-Reduction Methods: $\protect\theta %
_{1}=\max \{\protect\beta _{1},\protect\beta _{2}\}$.}
\begin{center}
\makebox{
\includegraphics[origin=bl,scale=.60,angle=0]{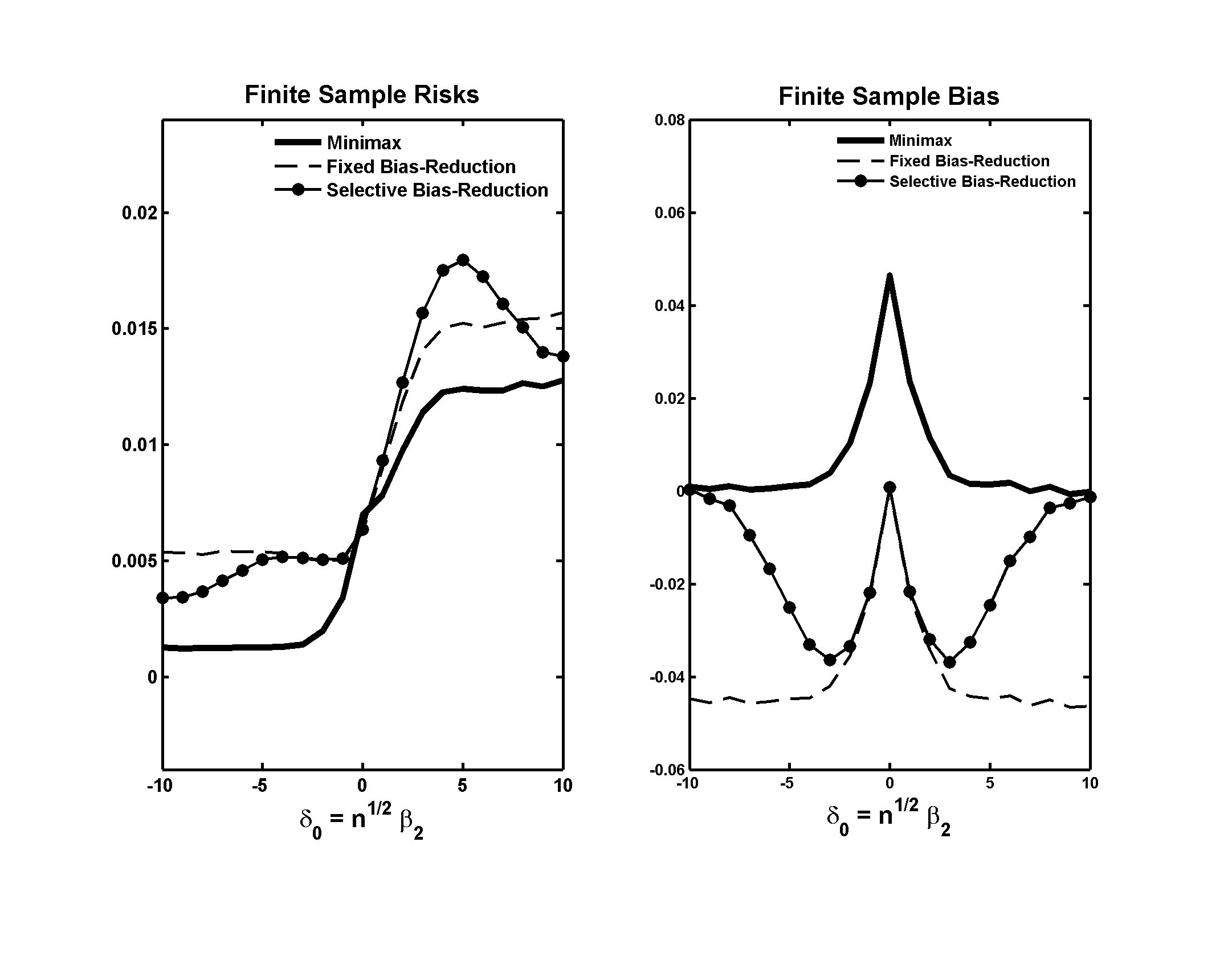}
}
\end{center}
\par
\label{figure}
\end{figure}

\subsection{Minimax Decision and Bias Reduction}

In the case of $\theta _{1}\equiv \max \{\beta _{1},\beta _{2}\}$, $%
b_{F}\equiv \mathbf{E}\left[ \max \{X_{11}-\beta _{1},X_{12}-\beta _{2}\}%
\right] $ becomes the asymptotic bias of the estimator $\hat{\theta}%
_{1}\equiv \max \{\hat{\beta}_{1},\hat{\beta}_{2}\}$ when $\beta _{1}=\beta
_{2}$. One may consider the following estimator of $b_{F}:$%
\begin{equation*}
\hat{b}_{F}\equiv \frac{1}{L}\sum_{i=1}^{L}\max \left( \hat{\Sigma}^{1/2}%
\mathbf{\xi }_{i}\right) \text{,}
\end{equation*}%
where $\mathbf{\xi }_{i}$ is drawn i.i.d. from $N(0,I_{2})$. This adjustment
term $\hat{b}_{F}$ is fixed over different values of $\beta _{2}-\beta _{1}$
(in large samples). Since the bias of $\max \{\hat{\beta}_{1},\hat{\beta}%
_{2}\}$ becomes prominent only when $\beta _{1}$ is close to $\beta _{2}$,
one may instead consider performing bias adjustment only when the estimated
difference $|\beta _{2}-\beta _{1}|$ is close to zero. Thus we also consider
the following estimated adjustment term:%
\begin{equation*}
\hat{b}_{S}\equiv \left( \frac{1}{L}\sum_{i=1}^{L}\max \left( \hat{\Sigma}%
^{1/2}\mathbf{\xi }_{i}\right) \right) 1\left\{ |\hat{\beta}_{2}-\hat{\beta}%
_{1}|<1.7/n^{1/3}\right\} .
\end{equation*}%
We compare the following two estimators with the minimax decision $\tilde{%
\theta}_{mx}$:%
\begin{equation*}
\hat{\theta}_{F}\equiv \max \{\hat{\beta}_{1},\hat{\beta}_{2}\}-\hat{b}_{F}/%
\sqrt{n}\text{ and }\hat{\theta}_{S}\equiv \max \{\hat{\beta}_{1},\hat{\beta}%
_{2}\}-\hat{b}_{S}/\sqrt{n}.
\end{equation*}%
We call $\hat{\theta}_{F}$ the estimator with fixed bias-reduction and $\hat{%
\theta}_{S}$ the estimator with selective bias-reduction. The results are
reported in Figure 1.

The finite sample risks of $\hat{\theta}_{F}$ are better than the minimax
decision $\hat{\theta}_{mx}$ only locally around $\delta _{0}=0$. The bias
reduction using $\hat{b}_{F}$ improves the estimator's performance in this
case. However, for other values of $\delta _{0}$, the bias reduction does
more harm than good because it lowers the bias when it is better not to, due
to increased variance. This is seen in the right-hand panel of Figure 1
which presents the finite sample bias of the estimators. With $\delta _{0}$
close to zero, the estimator with fixed bias-reduction eliminates the bias
almost entirely. However, for other values of $\delta _{0}$, this bias
correction induces negative bias, deteriorating the risk performances.

The estimator $\hat{\theta}_{S}$ with selective bias-reduction is designed
to be hybrid between the two extremes of $\hat{\theta}_{F}$ and $\tilde{%
\theta}_{mx}.$ When $\beta _{2}-\beta _{1}$ is estimated to be close to
zero, the estimator performs like $\hat{\theta}_{F}$ and when it is away
from zero, it performs like $\max \{\hat{\beta}_{1},\hat{\beta}_{2}\}$. As
expected, the bias of the estimator $\hat{\theta}_{S}$ is better than that
of $\hat{\theta}_{F}$ while successfully eliminating nearly the entire bias
when $\delta _{0}$ is close to zero. Nevertheless, it is remarkable that the
estimator shows highly unstable finite sample risk properties overall as
shown on the left panel in Figure 1. When $\delta _{0}$ is away from zero
and around 3 to 7, the performance is worse than the other estimators. This
result illuminates the fact that a reduction of bias does not always imply a
better risk performance.

The minimax decision shows finite sample risks that are robust over the
values of $\delta _{0}$. In fact, the estimated bias adjustment term $\hat{c}%
_{M_{1}}$ of the minimax decision is close to zero. This means that the
estimator $\hat{\theta}_{mx}$ requires zero bias adjustment, due to the
concern for its robust performance. In terms of finite sample bias, the
minimax estimator suffers from a substantially positive bias as compared to
the other two estimators, when $\delta _{0}$ is close to zero. The minimax
decision tolerates this bias because by doing so, it can maintain robust
performance for other cases where bias reduction is not needed. The minimax
estimator is ultimately concerned with the overall risk properties, not just
a bias component of the estimator, and as the left-hand panel of Figure 1
shows, it performs better than the other two estimators except when $\delta
_{0}$ is locally around zero, or when $\beta _{2}-\beta _{1}$ is around
roughly between $-0.057$ and $0.041$.

\section{Conclusion}

The paper proposes local asymptotic minimax estimators for a class of
nonregular parameters that are constructed by applying translation-scale
equivariant transform to a regular parameter. The results are extended to
the case where the nonregular parameters are transformed further by a
piecewise linear map with a single kink. The local asymptotic minimax
estimators take the form of a plug-in estimator with an additive bias
adjustment term. The bias adjustment term can be computed by a simulation
method. A small scale Monte Carlo simulation study demonstrates the robust
finite sample risk properties of the local asymptotic minimax estimators, as
compared to estimators based on alternative bias correction methods.

\section{Appendix: Mathematical Proofs}

\noindent \textsc{Proof of Lemma 1:} Property (a) follows immediately
because $g(\mathbf{0})=0$ by scale equivariance of $g$. Properties (b) and
(c) are due to translation equivariance of $g$. The first equality in
property (d) is due to scale equivariance of $g$, and the second equality
comes from the definition of directional derivatives. Lipschitz continuity
of $\tilde{g}(\mathbf{x};\cdot )$ on $\mathbf{R}^{d}$ stems from Lipschitz
continuity of $g$ (e.g. see the proof of Proposition 1.1 of Clarke (1998)).
Also, Lipschitz continuity of $g$ implies the uniform convergence on bounded
sets, because bounded directional differentiability and directional
differentiability in Assumption 1(i)(c) are equivalent when $g$ is a
Lipschitz map defined on a finite dimensional space. (See Shapiro (1990),
p.484.) $\blacksquare $\bigskip

\noindent \textsc{Proof of Lemma 2:} First, suppose to the contrary that $%
f_{1}(y)\neq f_{2}(y)$ for some $y\in \mathbf{R}$. Then since $f_{1}\circ
g_{1}=f_{2}\circ g_{2}$, it is necessary that $g_{1}(\mathbf{\beta })\neq
g_{2}(\mathbf{\beta })$ for some $\mathbf{\beta }\in \mathbf{R}^{d}$ such
that $g_{1}(\mathbf{\beta })=y,$ because $g_{1}(\mathbf{R}^{d})=\mathbf{R}$
and $g_{2}(\mathbf{R}^{d})=\mathbf{R}$, as we saw before. Hence%
\begin{equation}
(f_{1}\circ g_{1})(\mathbf{\beta })\neq (f_{2}\circ g_{1})(\mathbf{\beta }).
\label{noteq}
\end{equation}%
Now observe that $f_{2}(g_{1}(\mathbf{\beta }))=f_{2}(g_{2}(\mathbf{\beta }%
)+g_{1}(\mathbf{\beta })-g_{2}(\mathbf{\beta }))=f_{2}(g_{2}(\mathbf{\beta }%
+g_{1}(\mathbf{\beta })-g_{2}(\mathbf{\beta })))$. Since $f_{1}\circ
g_{1}=f_{2}\circ g_{2}$, the last term is equal to%
\begin{eqnarray*}
f_{1}(g_{1}(\mathbf{\beta }+g_{1}(\mathbf{\beta })-g_{2}(\mathbf{\beta })))
&=&f_{1}(2g_{1}(\mathbf{\beta })-g_{2}(\mathbf{\beta }))=f_{1}(g_{1}(2%
\mathbf{\beta }-g_{2}(\mathbf{\beta }))) \\
&=&f_{2}(g_{2}(2\mathbf{\beta }-g_{2}(\mathbf{\beta })))=f_{2}(g_{2}(\mathbf{%
\beta )})=f_{1}(g_{1}(\mathbf{\beta )}).
\end{eqnarray*}%
Therefore, we conclude that $f_{2}(g_{1}(\mathbf{\beta }))=f_{1}(g_{1}(%
\mathbf{\beta }))$ contradicting (\ref{noteq}).

Second, suppose to the contrary that $g_{1}(\mathbf{\beta })\neq g_{2}(%
\mathbf{\beta })$ for some $\mathbf{\beta }\in \mathbf{R}^{d}$ and $%
f_{1}=f_{2}$. First suppose that $g_{1}(\mathbf{\beta })>g_{2}(\mathbf{\beta 
})$. Fix arbitrary $a\in \mathbf{R}$ and $c\geq 0$ and let $c_{\Delta
}=c/\Delta _{1,2}(\mathbf{\beta })$ and $\Delta _{1,2}(\mathbf{\beta }%
)=g_{1}(\mathbf{\beta })-g_{2}(\mathbf{\beta })$. Then%
\begin{eqnarray*}
f_{1}(a+c) &=&f_{1}(a+\Delta _{1,2}(c_{\Delta }\mathbf{\beta }%
))=f_{1}(a+g_{2}(c_{\Delta }\mathbf{\beta })+\Delta _{1,2}(c_{\Delta }%
\mathbf{\beta })-g_{2}(c_{\Delta }\mathbf{\beta })) \\
&=&f_{1}(g_{2}(a+c_{\Delta }\mathbf{\beta }+\Delta _{1,2}(c_{\Delta }\mathbf{%
\beta })-g_{2}(c_{\Delta }\mathbf{\beta }))) \\
&=&f_{2}(g_{2}(a+c_{\Delta }\mathbf{\beta }+\Delta _{1,2}(c_{\Delta }\mathbf{%
\beta })-g_{2}(c_{\Delta }\mathbf{\beta }))) \\
&=&f_{1}(g_{1}(a+c_{\Delta }\mathbf{\beta }+\Delta _{1,2}(c_{\Delta }\mathbf{%
\beta })-g_{2}(c_{\Delta }\mathbf{\beta }))) \\
&=&f_{1}(a+g_{1}(c_{\Delta }\mathbf{\beta }-g_{2}(c_{\Delta }\mathbf{\beta }%
)+c_{\Delta }\Delta _{1,2}(\mathbf{\beta })\mathbf{)})=f_{1}(a+2c).
\end{eqnarray*}%
The choice of $a\in \mathbf{R}$ and $c\geq 0$ are arbitrary, and hence $%
f_{1}(\cdot )$ is constant on $\mathbf{R}$, contradicting the nonconstancy
condition for $f_{1}$.

Second, suppose that $g_{1}(\mathbf{\beta })<g_{2}(\mathbf{\beta })$. Then,
fix arbitrary $a\in \mathbf{R}$ and $c\leq 0$ and let $c_{\Delta }=c/\Delta
_{1,2}(\mathbf{\beta })$. Then similarly as before, we have%
\begin{eqnarray*}
f_{1}(a+c) &=&f_{1}(a+\Delta _{1,2}(c_{\Delta }\mathbf{\beta })) \\
&=&f_{1}(a+g_{1}(c_{\Delta }\mathbf{\beta }-g_{2}(c_{\Delta }\mathbf{\beta }%
)+c_{\Delta }\Delta _{1,2}(\mathbf{\beta })\mathbf{)})=f_{1}(a+2c),
\end{eqnarray*}%
because $\Delta _{1,2}(c_{\Delta }\mathbf{\beta })=c$. Therefore, again, $%
f_{1}(\cdot )$ is constant on $\mathbf{R}$, contradicting the nonconstancy
condition for $f_{1}$. $\blacksquare $\bigskip

We view convergence in distribution $\overset{d}{\rightarrow }$\ in the
proofs as convergence in $\mathbf{\bar{R}}^{d}$, so that the limit
distribution is allowed to be deficient in general. Choose $%
\{h_{i}\}_{i=1}^{m}$ from a complete orthonormal basis $\{h_{i}\}_{i=1}^{%
\infty }$ of $\bar{H}$. For $\mathbf{p}\in \mathbf{R}^{m}$, we consider $h(%
\mathbf{p})\equiv \Sigma _{i=1}^{m}p_{i}h_{i},$ $h_{i}\in H$, so that$\ \dot{%
\beta}_{j}(h(\mathbf{p}))=\sum_{i=1}^{m}\dot{\beta}_{j}(h_{i})p_{i},$ where $%
\dot{\beta}_{j}$ is the $j$-th element of $\mathbf{\dot{\beta}}.$ Let $%
\mathbf{B}$ be an $m\times d$ matrix such that%
\begin{equation}
\mathbf{B}\equiv \left[ 
\begin{tabular}{cccc}
$\dot{\beta}_{1}(h_{1})$ & $\dot{\beta}_{2}(h_{1})$ & $\cdots $ & $\dot{\beta%
}_{d}(h_{1})$ \\ 
$\dot{\beta}_{1}(h_{2})$ & $\dot{\beta}_{2}(h_{2})$ & $\cdots $ & $\dot{\beta%
}_{d}(h_{2})$ \\ 
$\vdots $ & $\vdots $ &  & $\vdots $ \\ 
$\dot{\beta}_{1}(h_{m})$ & $\dot{\beta}_{2}(h_{m})$ & $\cdots $ & $\dot{\beta%
}_{d}(h_{m})$%
\end{tabular}%
\right] .  \label{thetadot}
\end{equation}%
We assume that $m\geq d$ and $\mathbf{B}$ is a full column rank matrix.

We fix $h^{\prime }\in H$, and define $a_{i}=\langle h_{i},h^{\prime
}\rangle $ and $\mathbf{a}\in \mathbf{R}^{m}$ to be a column vector whose $i$%
-th entry is given by $a_{i}$. We also define $\mathbf{\zeta }\equiv (\zeta
(h_{1}),\cdot \cdot \cdot ,\zeta (h_{m}))^{\prime }$, where $\zeta $ is the
Gaussian process that appears in LAN, and with a small $\lambda >0$, let $%
F_{\lambda }(\cdot )$ be the cdf of $N(0,(I_{m}-\mathbf{a}(\mathbf{a}^{\top }%
\mathbf{a})^{-1}\mathbf{a}^{\top })/\lambda )$. Then by design, the
distribution of $h(\mathbf{p})\in \mathbf{R}$, with $\mathbf{p}\sim
F_{\lambda }$ concentrate on $\{h(\mathbf{p})\in \mathbf{R}:\langle h(%
\mathbf{p}),h^{\prime }\rangle =0\}$. Let $Z_{\lambda ,m}\in \mathbf{R}^{d}$
be a random vector following $N(0,\mathbf{B}^{\top }(I_{m}+\lambda (I_{m}-%
\mathbf{a}(\mathbf{\mathbf{a}^{\top }\mathbf{a}})^{-1}\mathbf{a}^{\top
})^{-1})^{-1}\mathbf{B})$.

Suppose that $\hat{\theta}\in \mathbf{R}$ is a sequence of estimators such
that along $\{P_{n,h^{\prime }}\}_{n\geq 1}$, with $h\subset H$ such that $%
\langle h,h^{\prime }\rangle =0$,%
\begin{equation*}
\left[ 
\begin{array}{c}
\sqrt{n}\{\hat{\theta}-g(\mathbf{\beta }_{n}(h+h^{\prime }))\} \\ 
\log dP_{n,h+h^{\prime }}/dP_{n,h^{\prime }}%
\end{array}%
\right] \overset{d}{\rightarrow }\left[ 
\begin{array}{c}
V-\tilde{g}_{0}(\mathbf{\dot{\beta}}(h)+\mathbf{r}))+\tilde{g}_{0}(\mathbf{r}%
)) \\ 
\zeta (h)-\frac{1}{2}\langle h,h\rangle%
\end{array}%
\right] ,
\end{equation*}%
for some nonstochastic vector $\mathbf{r}\in \mathbf{R}^{d}$, where $V\in 
\mathbf{R}$ is a random variable having a potentially deficient distribution
independent of $h\in H$.\footnote{%
Song (2014) on page 146 mistakenly refers to $V$ as a "random vector" in $%
\mathbf{R}^{d}$ when it is a random variable in $\mathbf{R}$. A similar
mistaken reference is found after the second display on page 149 of Song
(2014).} Let $\mathcal{L}_{g}^{h+h^{\prime }}$ be the limiting (potentially
deficient) distribution of $\sqrt{n}\{\hat{\theta}-g(\mathbf{\beta }%
_{n}(h+h^{\prime }))\}$ in $\mathbf{R}^{d}$ along $\{P_{n,h+h^{\prime
}}\}_{n\geq 1}$ for each $h\in H$ and$\ h^{\prime }\subset H$. The following
lemma is an adaptation of the generalized convolution theorem in van der
Vaart (1989).\bigskip

\noindent \textsc{Lemma\ A1:} \textit{Suppose that the map }$g$ \textit{%
satisfies Assumption 1(i) holds. Then the following holds.}

\noindent (i) \textit{For any} $\lambda >0,$ \textit{the distribution} $\int 
\mathcal{L}_{g}^{h(\mathbf{p})+h^{\prime }}dF_{\lambda }(\mathbf{p})$ 
\textit{is equal to that of} $-\tilde{g}_{0}(Z_{\lambda ,m}+W_{\lambda ,m}+%
\mathbf{r})+\tilde{g}_{0}(\mathbf{r})\in \mathbf{R},$ \textit{where\ }$%
W_{\lambda ,m}\in \mathbf{R}$\textit{\ is a random variable having a
potentially deficient distribution independent of }$Z_{\lambda ,m}$\textit{.}

\noindent (ii)\textit{\ As }$\lambda \rightarrow 0$\textit{\ first and then }%
$m\rightarrow \infty $\textit{, we have}$\ Z_{\lambda ,m}\overset{d}{%
\rightarrow }N(0,\Sigma )$.\bigskip

\noindent \textsc{Proof:}$\ $(i) Using Assumption 1(i) and applying Le Cam's
third lemma (van der Vaart and Wellner (1996), p.404), we find that for all $%
C\in \mathcal{B}(\mathbf{R}),$ the Borel $\sigma $-field of $\mathbf{R}$,%
\begin{eqnarray*}
\mathcal{L}_{g}^{h(\mathbf{p})+h^{\prime }}(C) &=&\mathbf{E}\left[ 1_{C}(V-%
\tilde{g}_{0}(\mathbf{B}^{\top }\mathbf{p}+\mathbf{r})+\tilde{g}_{0}(\mathbf{%
r}))e^{\mathbf{p}^{\top }\mathbf{\zeta }-\frac{1}{2}||\mathbf{p}||^{2}}%
\right] \\
&=&\mathbf{E}\left[ 1_{(-\tilde{g}_{0})^{-1}(C)}(-V+\mathbf{B}^{\top }%
\mathbf{p}+\mathbf{r-}\tilde{g}_{0}(\mathbf{r}))e^{\mathbf{p}^{\top }\mathbf{%
\zeta }-\frac{1}{2}||\mathbf{p}||^{2}}\right] ,
\end{eqnarray*}%
where $(-\tilde{g}_{0})^{-1}(C)\equiv \{\mathbf{x}\in \mathbf{R}^{d}:-\tilde{%
g}_{0}(\mathbf{x})\in C\}$. The second equality uses translation
equivariance of $\tilde{g}_{0}$. (See Lemma 1(c).) Define 
\begin{equation*}
\Sigma _{\lambda }\equiv \left( I_{m}+\lambda (I_{m}-\mathbf{a(\mathbf{%
\mathbf{a}^{\top }\mathbf{a}}})^{-1}\mathbf{a}^{\top })^{-1}\right) ^{-1}.
\end{equation*}%
Let $N_{\lambda }:\mathbf{R}^{m}\rightarrow \lbrack 0,1]$ be the
distribution function of $N(0,\Sigma _{\lambda })$. From the definition of $%
F_{\lambda }$, we write%
\begin{eqnarray*}
\int \mathcal{L}_{g}^{h(\mathbf{p})+h^{\prime }}(C)dF_{\lambda }(\mathbf{p})
&=&(2\pi )^{-m/2}\det (\lambda (I_{m}-\mathbf{a\mathbf{\mathbf{(\mathbf{%
\mathbf{a}^{\top }\mathbf{a}}}}})^{-1}\mathbf{a}^{\top }))^{-1/2} \\
&&\times \int \mathbf{E}\left[ 1_{(-\tilde{g}_{0})^{-1}(C)}\left( -V+\mathbf{%
B}^{\top }\mathbf{p}+\mathbf{r-}\tilde{g}_{0}(\mathbf{r})\right) e^{\mathbf{p%
}^{\top }\mathbf{\zeta }-\frac{1}{2}\mathbf{p}^{\top }\Sigma _{\lambda }^{-1}%
\mathbf{p}}\right] d\mathbf{p.}
\end{eqnarray*}%
By rearranging the terms and applying change of variables, we can rewrite
the integral as 
\begin{eqnarray*}
&&\int \mathbf{E}\left[ 1_{(-\tilde{g}_{0})^{-1}(C)}\left( -V+\mathbf{B}%
^{\top }\mathbf{p}+\mathbf{r-}\tilde{g}_{0}(\mathbf{r})\right) e^{-\frac{1}{2%
}(\mathbf{p-}\Sigma _{\lambda }\mathbf{\zeta })^{\top }\Sigma _{\lambda
}^{-1}(\mathbf{p-}\Sigma _{\lambda }\mathbf{\zeta })+\frac{1}{2}\mathbf{%
\zeta }^{\top }\Sigma _{\lambda }\mathbf{\zeta }}\right] d\mathbf{p} \\
&=&\int \mathbf{E}\left[ 1_{(-\tilde{g}_{0})^{-1}(C)}\left( -V+\mathbf{B}%
^{\top }\left( \mathbf{p}+\Sigma _{\lambda }\mathbf{\zeta }\right) +\mathbf{%
r-}\tilde{g}_{0}(\mathbf{r})\right) e^{-\frac{1}{2}\mathbf{p}\Sigma
_{\lambda }^{-1}\mathbf{p}+\frac{1}{2}\mathbf{\zeta }^{\top }\Sigma
_{\lambda }\mathbf{\zeta }}\right] d\mathbf{p.}
\end{eqnarray*}%
Therefore, we conclude that%
\begin{equation*}
\int \mathcal{L}_{g}^{h(\mathbf{p})+h^{\prime }}(C)dF_{\lambda }(\mathbf{p}%
)=\int \mathbf{E}\left[ 1_{(-\tilde{g}_{0})^{-1}(C)}\left( -V+\mathbf{B}%
^{\top }\left( \mathbf{p}+\Sigma _{\lambda }\mathbf{\zeta }\right) +\mathbf{%
r-}\tilde{g}_{0}(\mathbf{r})\right) c_{\lambda }(\mathbf{\zeta })\right]
dN_{\lambda }(\mathbf{p),}
\end{equation*}%
where $c_{\lambda }(\mathbf{\zeta })\equiv e^{\frac{1}{2}\mathbf{\zeta }%
^{\top }\Sigma _{\lambda }\mathbf{\zeta }}\cdot \det (\lambda (I_{m}-\mathbf{%
a\mathbf{(\mathbf{\mathbf{a}^{\top }\mathbf{a}}}})^{-1}\mathbf{a}^{\top
}))^{-1/2}/\det (\Sigma _{\lambda })^{-1/2}.$ When we let $W_{\lambda ,m}$
be a random variable having potentially deficient distribution $\mathcal{W}%
_{\lambda ,m}$ defined by%
\begin{equation*}
\mathcal{W}_{\lambda ,m}(C)\equiv \mathbf{E}\left[ 1_{(-\tilde{g}%
_{0})^{-1}(C)}\left( V-\mathbf{B}^{\top }\Sigma _{\lambda }\mathbf{\zeta }%
\right) c_{\lambda }(\mathbf{\zeta })\right] ,\ C\in \mathcal{B}(\mathbf{R}),
\end{equation*}%
the distribution $\int \mathcal{L}_{g}^{h(\mathbf{p})+h^{\prime
}}dF_{\lambda }(\mathbf{p})$ is equal to that of $-\tilde{g}_{0}(Z_{\lambda
,m}+W_{\lambda ,m}+\mathbf{r}-\tilde{g}_{0}(\mathbf{r}))$.

\noindent (ii) Since the sequence $\{h_{i}\}_{i=1}^{\infty }$ is a complete
orthonormal basis of $\bar{H}$, the covariance matrix of $Z_{\lambda ,m}$
converges to $\Sigma $ as $\lambda \rightarrow 0$ and then $m\rightarrow
\infty $. $\blacksquare $\bigskip

We introduce some notation. Define $||\cdot ||_{BL}$ on the space of Borel
measurable functions on $\mathbf{R}^{d}:$%
\begin{equation*}
||f||_{BL}\equiv \sup_{x\neq y}|f(x)-f(y)|/||x-y||+\sup_{x}|f(x)|.
\end{equation*}%
For any two probability measures $P$ and $Q$ on $\mathcal{B}(\mathbf{R}%
^{d}), $ define 
\begin{equation}
d_{\mathcal{P}}(P,Q)\equiv \sup \left\{ \left\vert \int fdP-\int
fdQ\right\vert :||f||_{BL}\leq 1\right\} .  \label{metric}
\end{equation}%

\noindent \textsc{Lemma A4 :} \textit{Suppose that for each }$n\geq 1,$%
\textit{\ }$\{P_{n,h}:h\in H\}$\textit{\ is the set of probability measures
indexed by a Hilbert space }$(H,\langle \cdot ,\cdot \rangle )$, \textit{%
such that for each }$h\in H,$%
\begin{equation*}
\log \frac{dP_{n,h}}{dP_{n,0}}=\zeta _{n}(h)-\frac{1}{2}\langle h,h\rangle ,%
\text{ \textit{under} }\{P_{n,0}\},
\end{equation*}%
\textit{where for each }$h,h^{\prime }\in H$, $[\zeta _{n}(h),\zeta
_{n}(h^{\prime })]\overset{d}{\rightarrow }[\zeta (h),\zeta (h^{\prime })]$ 
\textit{under }$\{P_{n,0}\}$\textit{, and }$\zeta (\cdot )$ \textit{is a
Gaussian process on }$H$ \textit{with covariance function }$\mathbf{E}[\zeta
(h_{1})\zeta (h_{2})]=\langle h_{1},h_{2}\rangle $, $h_{1},h_{2}\in H$.

\textit{Then for each }$h,h^{\prime }\in H$\textit{\ such that }$\langle
h,h^{\prime }\rangle =0$,%
\begin{equation}
\log \frac{dP_{n,h+h^{\prime }}}{dP_{n,h^{\prime }}}\overset{d}{\rightarrow }%
\ \ \zeta (h)-\frac{1}{2}\langle h,h\rangle ,\text{ \textit{under} }%
\{P_{n,h^{\prime }}\}.  \label{dc}
\end{equation}%
\noindent \textsc{Proof :} Since 
\begin{equation*}
\log dP_{n,h+h^{\prime }}/dP_{n,h^{\prime }}=\log dP_{n,h+h^{\prime
}}/dP_{n,0}-\log dP_{n,h^{\prime }}/dP_{n,0},
\end{equation*}%
we observe that by the condition of the lemma,%
\begin{equation*}
\left[ 
\begin{array}{c}
\log dP_{n,h+h^{\prime }}/dP_{n,h^{\prime }} \\ 
\log dP_{n,h^{\prime }}/dP_{n,0}%
\end{array}%
\right] \overset{d}{\rightarrow }\left[ 
\begin{array}{c}
\zeta (h+h^{\prime })-\zeta (h^{\prime })-\frac{1}{2}\langle h,h\rangle \\ 
\zeta (h^{\prime })-\frac{1}{2}\langle h^{\prime },h^{\prime }\rangle%
\end{array}%
\right] ,
\end{equation*}%
under $\{P_{n,0}\}_{n\geq 1}$, because $\langle h,h^{\prime }\rangle =0$. By
Le Cam's third lemma (van der Vaart and Wellner (1996), p.404), under$%
\{P_{n,h^{\prime }}\}_{n\geq 1}$,%
\begin{equation*}
\log dP_{n,h+h^{\prime }}/dP_{n,h^{\prime }}\overset{d}{\rightarrow }%
\mathcal{L}\text{,}
\end{equation*}%
where $\mathcal{L}$ is a probability measure on $\mathcal{B}(\mathbf{R})$
such that for any $B\in \mathcal{B}(\mathbf{R}),$%
\begin{eqnarray*}
\mathcal{L}(B) &=&\mathbf{E}\left[ 1_{B}\left( \zeta (h+h^{\prime })-\zeta
(h^{\prime })-\frac{1}{2}\langle h,h\rangle \right) e^{\zeta (h^{\prime })-%
\frac{1}{2}\langle h^{\prime },h^{\prime }\rangle }\right] \\
&=&\mathbf{E}\left[ 1_{B}\left( \zeta (h+h^{\prime })-\zeta (h^{\prime })-%
\frac{1}{2}\langle h,h\rangle \right) \right] \\
&=&\mathbf{E}\left[ 1_{B}\left( \zeta (h)-\frac{1}{2}\langle h,h\rangle
\right) \right] .
\end{eqnarray*}%
The second equality above follows because $\zeta (h+h^{\prime })-\zeta
(h^{\prime })$ and $\zeta (h^{\prime })$ are independent and $\mathbf{E}%
[e^{\zeta (h^{\prime })-\frac{1}{2}\langle h^{\prime },h^{\prime }\rangle
}]=1$, and the third equality above follows because 
\begin{equation*}
\zeta (h+h^{\prime })-\zeta (h^{\prime })\overset{d}{=}N(0,||h||^{2})\overset%
{d}{=}\zeta (h).
\end{equation*}%
Hence we obtain (\ref{dc}). $\blacksquare $\bigskip

\noindent \textsc{Proof of Lemma 3:} We show that%
\begin{eqnarray}
&&\sup_{b\in \lbrack 0,\infty )}\ \underset{n\rightarrow \infty }{\text{%
liminf}}\sup_{h\in H_{n,b}}\mathbf{E}_{h}\left[ \tau (|\sqrt{n}\{\hat{\theta}%
-g(\mathbf{\beta }_{n}(h))\}|)\right]  \label{step} \\
&\geq &\sup_{\mathbf{r}\in \Gamma }\int \mathbf{E}\left[ \tau (|g\left( Z+%
\mathbf{r}\right) +w|)\right] dF(w),  \notag
\end{eqnarray}%
for some $F\in \mathcal{F}$. Then the proof is complete by taking infimum
over $F\in \mathcal{F}$.

First, we choose $\mathbf{r}\in \mathbf{R}^{d}$. Then we can find some $%
h^{\prime }\in \overline{H}$ such that $\mathbf{r=\dot{\beta}(}h^{\prime }%
\mathbf{)}$. More specifically, let $\mathbf{q}=\Sigma ^{-1}\mathbf{r}$ and
define $h^{\prime }=\sum_{i=1}^{d}\tilde{\beta}_{i}q_{i}$, where for each $%
i=1,\cdot \cdot \cdot ,d$, $\tilde{\beta}_{i}\in \overline{H}$ is such that $%
\langle \tilde{\beta}_{i},h\rangle =\mathbf{e}_{i}^{\top }\mathbf{\dot{\beta}%
(}h\mathbf{)}$ for all $h\in H$, and $q_{i}$ is the $i$-th entry of $\mathbf{%
q}$. Then for this choice of $h^{\prime }$, we can show that $\mathbf{r=\dot{%
\beta}(}h^{\prime }\mathbf{)}$.

Fix $b/2\geq ||h^{\prime }||\cdot ||\mathbf{\dot{\beta}}^{\ast }||$, where $%
\mathbf{\dot{\beta}}^{\ast }=(\dot{\beta}_{\mathbf{e}_{1}}^{\ast },\cdot
\cdot \cdot ,\dot{\beta}_{\mathbf{e}_{d}}^{\ast })^{\top }$ and $\dot{\beta}%
_{\mathbf{e}_{m}}^{\ast }$'s are as defined after Assumption 2. We note that%
\begin{eqnarray}
&&\underset{n\rightarrow \infty }{\text{liminf}}\sup_{h\in H_{n,b}}\mathbf{E}%
_{h}\left[ \tau (|\sqrt{n}\{\hat{\theta}-g(\mathbf{\beta }_{n}(h))\}|)\right]
\label{ineq4} \\
&\geq &\ \underset{n\rightarrow \infty }{\text{liminf}}\sup_{h\in
H_{n,b/2}^{\ast }}\mathbf{E}_{h+h^{\prime }}\left[ \tau (|\sqrt{n}\{\hat{%
\theta}-g(\mathbf{\beta }_{n}(h+h^{\prime }))\}|)\right] ,  \notag
\end{eqnarray}%
where 
\begin{equation*}
H_{n,b/2}^{\ast }\equiv \left\{ h\in H_{n,b/2}:\langle h,h^{\prime }\rangle
=0\right\} .
\end{equation*}%
As in the proof of Theorem 3.11.5 of van der Vaart and Wellner (1996)
(p.417), choose an orthonormal basis $\{h_{i}\}_{i=1}^{\infty }$ from $\bar{H%
}$. We fix $m$ and take $\{h_{i}\}_{i=1}^{m}\subset H$ and consider $h(%
\mathbf{p})=\sum p_{i}h_{i}$ for some $\mathbf{p}=(p_{i})_{i=1}^{m}\in 
\mathbf{R}^{m}\ $such that $h(\mathbf{p})\in H.$ Fix $\lambda >0$ and let $%
F_{\lambda }(\mathbf{p})$ be as defined prior to Lemma A1 (with $h^{\prime
}\in H$ chosen previously in this proof.) Note that by design, any vector $%
\mathbf{p}$ in the support of the distribution $F_{\lambda }$ satisfies that 
$\langle h(\mathbf{p}),h^{\prime }\rangle =0$. Hence note that for fixed $%
M>0,$%
\begin{eqnarray}
&&\underset{n\rightarrow \infty }{\text{liminf}}\sup_{h\in H_{n,b/2}^{\ast }}%
\mathbf{E}_{h+h^{\prime }}\left[ \tau (|\sqrt{n}\{\hat{\theta}-g(\mathbf{%
\beta }_{n}(h+h^{\prime }))\}|)\right]  \label{inte} \\
&\geq &\underset{n\rightarrow \infty }{\ \text{liminf}}\int \mathbf{E}_{h(%
\mathbf{p})+h^{\prime }}\left[ \tau _{M}\left( \left\vert V_{n,h(\mathbf{p}%
)+h^{\prime }}\right\vert \right) \right] 1\left\{ h(\mathbf{p})\in
H_{n,b/2}^{\ast }\right\} dF_{\lambda }(\mathbf{p})  \notag \\
&\geq &\int \underset{n\rightarrow \infty }{\text{liminf}}\ \mathbf{E}_{h(%
\mathbf{p})+h^{\prime }}\left[ \tau _{M}\left( \left\vert V_{n,h(\mathbf{p}%
)+h^{\prime }}\right\vert \right) \right] dF_{\lambda }(\mathbf{p})  \notag
\\
&&-M\underset{n\rightarrow \infty }{\text{limsup}}\int 1\left\{ h(\mathbf{p}%
)\notin H_{n,b/2}^{\ast }\right\} dF_{\lambda }(\mathbf{p}),  \notag
\end{eqnarray}%
where $V_{n,h}\equiv \sqrt{n}\{\hat{\theta}-g(\mathbf{\beta }_{n}(h))\}$.
The second inequality uses Fatou's lemma.

We write%
\begin{eqnarray*}
&&\sqrt{n}\{\hat{\theta}-g(\mathbf{\beta }_{n}(h+h^{\prime }))\} \\
&=&\sqrt{n}\{\hat{\theta}-g(\mathbf{\beta }_{n}(h^{\prime }))\}-\sqrt{n}g(%
\mathbf{\beta }_{n}(h+h^{\prime }))+\sqrt{n}g(\mathbf{\beta }_{n^{\prime
}}(h^{\prime })).
\end{eqnarray*}%
Then%
\begin{eqnarray}
&&\sqrt{n}g(\mathbf{\beta }_{n}(h+h^{\prime }))-\sqrt{n}g(\mathbf{\beta }%
_{n}(0))  \label{der1} \\
&=&\sqrt{n}\left\{ g(\mathbf{\beta }_{n}(h+h^{\prime })-\mathbf{\beta }%
_{n}(0)+\mathbf{\beta }_{n}(0)))-g(\mathbf{\beta }_{n}(0))\right\}  \notag \\
&=&\sqrt{n}\{g(\mathbf{\dot{\beta}}(h+h^{\prime })/\sqrt{n}+\mathbf{\beta }%
_{n}(0)))-g(\mathbf{\beta }_{n}(0))\}+o(1)  \notag \\
&=&\tilde{g}(\mathbf{\beta }_{n}(0);\mathbf{\dot{\beta}}(h+h^{\prime
}))+o(1)=\tilde{g}(\mathbf{\beta }_{n}(0);\mathbf{\dot{\beta}}(h)+\mathbf{r}%
)+o(1)  \notag \\
&=&\tilde{g}_{0}(\mathbf{\dot{\beta}}(h)+\mathbf{r}))+o(1),  \notag
\end{eqnarray}%
where the second to the last equality follows by the linearity of $\mathbf{%
\dot{\beta}}$ and the choice of $\mathbf{r}$, and the last equality follows
because $\mathbf{\beta }_{n}(0)=\mathbf{\beta }(P_{\alpha _{0}})=\mathbf{%
\beta }_{0}$ and by the definition of $\tilde{g}_{0}(\cdot )$. Similarly,%
\begin{equation}
\sqrt{n}g(\mathbf{\beta }_{n^{\prime }}(h^{\prime }))-\sqrt{n}g(\mathbf{%
\beta }_{n}(0))=\tilde{g}_{0}(\mathbf{r})+o(1).  \label{der2}
\end{equation}%
Combining (\ref{der1}) and (\ref{der2}), we find that 
\begin{equation}
\sqrt{n}g(\mathbf{\beta }_{n}(h+h^{\prime }))-\sqrt{n}g(\mathbf{\beta }%
_{n}(h^{\prime }))\rightarrow \tilde{g}_{0}(\mathbf{\dot{\beta}}(h)+\mathbf{r%
})-\tilde{g}_{0}(\mathbf{r}),  \label{cv32}
\end{equation}%
as $n\rightarrow \infty $.

Applying Prohorov's Theorem (in $\mathbf{\bar{R}}^{d}$), we find that for
any subsequence of $\{n\}$, there exists a further subsequence $\{n^{\prime
}\}$ along which (under $\{P_{n^{\prime },h^{\prime }}\}$) 
\begin{equation*}
\sqrt{n^{\prime }}\{\hat{\theta}-g(\mathbf{\beta }_{n}(h^{\prime }))\}%
\overset{d}{\rightarrow }V,
\end{equation*}%
where $V\in \mathbf{\bar{R}}$ is a random variable having a potentially
deficient distribution. Observe that 
\begin{eqnarray}
\sqrt{n^{\prime }}\{\hat{\theta}-g(\mathbf{\beta }_{n^{\prime }}(h+h^{\prime
}))\} &=&\sqrt{n^{\prime }}\{\hat{\theta}-g(\mathbf{\beta }_{n^{\prime
}}(h^{\prime }))\}  \label{der4} \\
&&-\{\sqrt{n}g(\mathbf{\beta }_{n}(h+h^{\prime }))-\sqrt{n}g(\mathbf{\beta }%
_{n}(h^{\prime }))\}  \notag \\
&&\overset{d}{\rightarrow }V-\tilde{g}_{0}(\mathbf{\dot{\beta}}(h)+\mathbf{r}%
)+\tilde{g}_{0}(\mathbf{r}).  \notag
\end{eqnarray}%
Invoking Assumption 2, Lemma A4, and (\ref{der4}), and noting that marginal
tightness implies joint tightness, we apply Prohorov's Theorem to deduce
that for any subsequence of $\{n\}$, there exists a further subsequence $%
\{n^{\prime }\}$ along which $\mathbf{r}_{n^{\prime }}\rightarrow \mathbf{%
r\equiv \beta }(h^{\prime })$, and (under $P_{n^{\prime },h^{\prime }}$) 
\begin{equation*}
\left[ 
\begin{array}{c}
\sqrt{n^{\prime }}\{\hat{\theta}-g(\mathbf{\beta }_{n^{\prime }}(h+h^{\prime
}))\} \\ 
\log dP_{n^{\prime },h+h^{\prime }}/dP_{n^{\prime },h^{\prime }}%
\end{array}%
\right] \overset{d}{\rightarrow }\left[ 
\begin{array}{c}
V-\tilde{g}_{0}(\mathbf{\dot{\beta}}(h)+\mathbf{r})+\tilde{g}_{0}(\mathbf{r})
\\ 
\zeta (h)-\frac{1}{2}\langle h,h\rangle%
\end{array}%
\right] ,
\end{equation*}%
where $\sqrt{n^{\prime }}\{\hat{\theta}-g(\mathbf{\beta }_{n^{\prime
}}(h^{\prime }))\}\overset{d}{\rightarrow }V$ under $P_{n^{\prime
},h^{\prime }}$. By Lemma A1,%
\begin{equation*}
\int \underset{n\rightarrow \infty }{\ \text{liminf}}\ \mathbf{E}_{h(\mathbf{%
p})+h^{\prime }}\left[ \tau _{M}\left( \left\vert V_{n,h(\mathbf{p}%
)+h^{\prime }}\right\vert \right) \right] dF_{\lambda }(\mathbf{p})=\mathbf{E%
}[\tau _{M}(|\tilde{g}_{0}(Z_{\lambda ,m}+W_{\lambda ,m}+\mathbf{r})-\tilde{g%
}_{0}(\mathbf{r})|)],
\end{equation*}%
where $Z_{\lambda ,m}$ is as defined prior to Lemma A1 and $W_{\lambda
,m}\in \mathbf{R}$ is a random variable having a potentially deficient
distribution and independent of $Z_{\lambda ,m}$. Furthermore, by Assumption
2 (regularity of $\mathbf{\beta }_{n}(h)$), we have for each $\mathbf{p}\in 
\mathbf{R}^{m}$,%
\begin{equation*}
1\left\{ h(\mathbf{p})\in H_{n,b/2}^{\ast }\right\} \rightarrow 1\left\{ h(%
\mathbf{p})\in H_{b/2}^{\ast }\right\} ,
\end{equation*}%
as $n\rightarrow \infty $, where $H_{b}^{\ast }\equiv \{h\in H:||\mathbf{%
\dot{\beta}}(h)||\leq b,\langle h,h^{\prime }\rangle =0\},$ and as $%
b\uparrow \infty $,%
\begin{equation*}
1\left\{ h(\mathbf{p})\in H_{b/2}^{\ast }\right\} \rightarrow 1\left\{ h(%
\mathbf{p})\in H^{\ast }\right\} ,
\end{equation*}%
where $H^{\ast }\equiv \{h\in H:\langle h,h^{\prime }\rangle =0\}$.
Therefore, since for each $\mathbf{p}$ in the support of $F_{\lambda }$, we
have $h(\mathbf{p})\in H^{\ast }$, we send $n\rightarrow \infty $ and $%
b\uparrow \infty ,$ and apply the Dominated Convergence Theorem to conclude
that%
\begin{equation*}
\underset{b\rightarrow \infty }{\ \text{lim}}\ \underset{n\rightarrow \infty 
}{\text{limsup}}\int 1\left\{ h(\mathbf{p})\notin H_{n,b/2}^{\ast }\right\}
dF_{\lambda }(\mathbf{p})=0\text{.}
\end{equation*}%
Thus, we conclude from (\ref{inte}) that%
\begin{eqnarray}
&&\underset{b\rightarrow \infty }{\ \text{lim}}\ \underset{n\rightarrow
\infty }{\text{liminf}}\sup_{h\in H_{n,b/2}^{\ast }}\mathbf{E}_{h+h^{\prime
}}\left[ \tau (|\sqrt{n}\{\hat{\theta}-g(\mathbf{\beta }_{n}(h+h^{\prime
}))\}|)\right]  \label{ineq3} \\
&\geq &\mathbf{E}[\tau _{M}(|\tilde{g}_{0}(Z_{\lambda ,m}+W_{\lambda ,m}+%
\mathbf{r})-\tilde{g}_{0}(\mathbf{r})|)].  \notag
\end{eqnarray}%
By Lemma A1(ii), as $\lambda \rightarrow 0$ and then $m\rightarrow \infty $, 
$Z_{\lambda ,m}$ converges in distribution to $Z$. Since $\{[Z_{\lambda
,m}^{\top },W_{\lambda ,m}]^{\top }\in \mathbf{\bar{R}}^{d+1}:(\lambda
,m)\in (0,\infty )\times \{1,2,\cdot \cdot \cdot \}\}$ is uniformly tight in 
$\mathbf{\bar{R}}^{d+1}$, by Prohorov's Theorem, for any subsequence of $%
\{\lambda _{k}\}_{k=1}^{\infty }$ with $\lambda _{k}\rightarrow 0$ as $%
k\rightarrow \infty $, and subsequence of $\{m\}$, there exist further
subsequences $\{\lambda _{k^{\prime }}\}\subset \{\lambda _{k}\}$ and $%
\{m^{\prime }\}\subset \{m\}$, such that as $k^{\prime }\rightarrow 0$ and
then $m^{\prime }\rightarrow \infty $, 
\begin{equation*}
\lbrack Z_{\lambda _{k^{\prime }},m}^{\top },W_{\lambda _{k^{\prime
}},m^{\prime }}]^{\top }\overset{d}{\rightarrow }[Z^{\top },W]^{\top },
\end{equation*}%
for some random variable $W_{m}$ having a potentially deficient
distribution. By applying this to the right hand side of (\ref{ineq3}) and
recalling (\ref{ineq4}), and noting that the choice of $\mathbf{r}\in 
\mathbf{R}^{d}$ was arbitrary, we conclude that%
\begin{eqnarray}
&&\underset{b\uparrow \infty }{\text{lim}}\ \underset{n\rightarrow \infty }{%
\text{liminf}}\sup_{h\in H_{n,b}}\mathbf{E}_{h}\left[ \tau (|\sqrt{n}\{\hat{%
\theta}-g(\mathbf{\beta }_{n}(h))\}|)\right]  \label{bound} \\
&\geq &\sup_{\mathbf{r}\in \mathbf{R}^{d}}\int \mathbf{E}\left[ \tau _{M}(|%
\tilde{g}_{0}(Z+w+\mathbf{r})-\tilde{g}_{0}(\mathbf{r})|)\right] dF(w), 
\notag
\end{eqnarray}%
where $F$ is an element of $\mathcal{F}^{\ast }$ and $\mathcal{F}^{\ast }$
is the collection of distributions on $\mathcal{B}(\mathbf{\bar{R}})$.

Fix $F\in \mathcal{F}^{\ast }$. As for the last integral in (\ref{bound}),
we write it as 
\begin{eqnarray}
&&\int \mathbf{E}\left[ \tau _{M}(|\tilde{g}_{0}(Z+\mathbf{r})-\tilde{g}_{0}(%
\mathbf{r})+w|)\right] dF(w)  \label{dec} \\
&=&\int \mathbf{E}\left[ \tau _{M}(|\tilde{g}_{0}(Z+\mathbf{r})-\tilde{g}%
_{0}(\mathbf{r})+w|)1\{w\in \mathbf{\bar{R}}\backslash \mathbf{R}\}\right]
dF(w)  \notag \\
&&+\int \mathbf{E}\left[ \tau _{M}(|\tilde{g}_{0}(Z+\mathbf{r})-\tilde{g}%
_{0}(\mathbf{r})+w|)1\{w\in \mathbf{R}\}\right] dF(w).  \notag
\end{eqnarray}%
Since $\tilde{g}_{0}(Z+\mathbf{r})-\tilde{g}_{0}(\mathbf{r})\in \mathbf{R}$,
for $w\in \mathbf{\bar{R}}\backslash \mathbf{R},$ 
\begin{equation*}
\tau _{M}\left( |\tilde{g}_{0}(Z+\mathbf{r})-\tilde{g}_{0}(\mathbf{r}%
)+w|\right) =\min \left\{ \sup_{x\in \lbrack 0,\infty )}\tau (x),M\right\} ,
\end{equation*}%
so that%
\begin{eqnarray*}
&&\int \mathbf{E}\left[ \tau _{M}\left( |\tilde{g}_{0}(Z+\mathbf{r})-\tilde{g%
}_{0}(\mathbf{r})+w|\right) 1\{w\in \mathbf{\bar{R}}\backslash \mathbf{R}\}%
\right] dF(w) \\
&=&\min \left\{ \sup_{x\in \lbrack 0,\infty )}\tau (x),M\right\} \cdot \int_{%
\mathbf{\bar{R}}\backslash \mathbf{R}}dF(w).
\end{eqnarray*}%
From (\ref{dec}), we conclude that 
\begin{eqnarray*}
&&\int \mathbf{E}\left[ \tau _{M}(|\tilde{g}_{0}(Z+\mathbf{r})-\tilde{g}_{0}(%
\mathbf{r})+w|)\right] dF(w) \\
&=&\min \left\{ \sup_{x\in \lbrack 0,\infty )}\tau (x),M\right\} \cdot \int_{%
\mathbf{\bar{R}}\backslash \mathbf{R}}dF(w) \\
&&+\int \mathbf{E}\left[ \tau _{M}(|\tilde{g}_{0}(Z+\mathbf{r})-\tilde{g}%
_{0}(\mathbf{r})+w|)1\{w\in \mathbf{R}\}\right] dF(w).
\end{eqnarray*}

We identify $\mathcal{F}$ as the subset of $\mathcal{F}^{\ast }$ such that
for each $F\in \mathcal{F}$, $\int_{\mathbf{\bar{R}}\backslash \mathbf{R}%
}dF(w)=0$ and $\int_{\mathbf{R}}dF(w)=1$. Since 
\begin{equation*}
\tau _{M}(|\tilde{g}_{0}(Z+\mathbf{r})-\tilde{g}_{0}(\mathbf{r})+w|)\leq
\min \left\{ \sup_{x\in \lbrack 0,\infty )}\tau (x),M\right\} \text{
everywhere,}
\end{equation*}%
the lower bound in (\ref{bound}) remains the same if we replace $\mathcal{F}%
^{\ast }$ by $\mathcal{F}$. Since $\tau _{M}$ increases in $M$, we obtain
the desired bound by sending $M\uparrow \infty $. $\blacksquare $\bigskip

For given $M_{1},a>0$ and $c\in \mathbf{R}$, define%
\begin{equation}
B_{M_{1}}(c;a)\equiv \sup_{\mathbf{r}\in \mathbf{R}^{d}}\mathbf{E}\left[
\tau _{M_{1}}(a|\tilde{g}_{0}(Z+\mathbf{r})-\tilde{g}_{0}(\mathbf{r})+c|)%
\right] ,  \label{BM}
\end{equation}%
and%
\begin{equation*}
E_{M_{1}}(a)\equiv \left\{ c\in \lbrack -M_{1},M_{1}]:B_{M_{1}}(c;a)\leq
\inf_{c_{1}\in \lbrack -M_{1},M_{1}]}B_{M_{1}}(c_{1};a)\right\} .
\end{equation*}%
Let $c_{M_{1}}^{\ast }(a)\equiv \sup E_{M_{1}}(a)$. We also define 
\begin{equation*}
\bar{g}_{n}(\mathbf{z})\equiv g\left( \mathbf{z+}\varepsilon _{n}^{-1}(%
\mathbf{\beta }_{0}-g(\mathbf{\beta }_{0}))\right) ,
\end{equation*}%
for $\mathbf{z\in R}^{d}$, and%
\begin{eqnarray*}
\bar{B}_{M_{1}}(c;a) &\equiv &\sup_{\mathbf{r}\in \lbrack -M_{1},M_{1}]^{d}}%
\frac{1}{L}\sum_{i=1}^{L}\tau _{M_{1}}\left( a\left\vert \bar{g}_{n}(\hat{%
\Sigma}^{1/2}\mathbf{\xi }_{i}+\mathbf{r})-\bar{g}_{n}(\mathbf{r}%
)+c\right\vert \right) \text{,} \\
\tilde{B}_{M_{1}}(c;a) &\equiv &\sup_{\mathbf{r}\in \lbrack
-M_{1},M_{1}]^{d}}\frac{1}{L}\sum_{i=1}^{L}\tau _{M_{1}}\left( a\left\vert 
\bar{g}_{n}(\Sigma ^{1/2}\mathbf{\xi }_{i}+\mathbf{r})-\bar{g}_{n}(\mathbf{r}%
)+c\right\vert \right) \text{,}
\end{eqnarray*}%
and%
\begin{equation*}
B_{M_{1}}^{\ast }(c;a)\equiv \sup_{\mathbf{r}\in \lbrack -M_{1},M_{1}]^{d}}%
\mathbf{E}\left[ \tau _{M_{1}}\left( a\left\vert \bar{g}_{n}(\Sigma ^{1/2}%
\mathbf{\xi }_{i}+\mathbf{r})-\bar{g}_{n}(\mathbf{r})+c\right\vert \right) %
\right] .
\end{equation*}%
We also define%
\begin{equation*}
E_{M_{1}}^{\ast }(a)\equiv \left\{ c\in \lbrack
-M_{1},M_{1}]:B_{M_{1}}^{\ast }(c;a)\leq \inf_{c_{1}\in \lbrack
-M_{1},M_{1}]}B_{M_{1}}^{\ast }(c_{1};a)\right\} .
\end{equation*}

\noindent \textsc{Lemma A5}: \textit{Suppose that Assumptions} 1(i), 4,%
\textit{\ and }5 \textit{hold. Then as }$M\rightarrow \infty ,$%
\begin{equation*}
\lim_{n\rightarrow \infty }\sup_{h\in H}P_{n,h}\left\{ \sup_{c\mathbf{\in }%
[-M_{1},M_{1}]}\left\vert B_{M_{1}}^{\ast }(c;a)-\hat{B}_{M_{1}}(c;a)\right%
\vert >M(L^{-1/2}+n^{-1/2}\varepsilon _{n}^{-1})\right\} \rightarrow 0.
\end{equation*}

\noindent \textsc{Proof}: Note that%
\begin{eqnarray*}
\left\vert \bar{g}_{n}(\mathbf{z})-\hat{g}_{n}(\mathbf{z})\right\vert
&=&\left\vert g\left( \mathbf{z+}\varepsilon _{n}^{-1}(\mathbf{\beta }_{0}-g(%
\mathbf{\beta }_{0}))\right) -g(\mathbf{z+}\varepsilon _{n}^{-1}(\mathbf{%
\hat{\beta}}-g(\mathbf{\hat{\beta}})))\right\vert \\
&\leq &2\varepsilon _{n}^{-1}\left\Vert \mathbf{\beta }_{0}-\mathbf{\hat{%
\beta}}\right\Vert ,
\end{eqnarray*}%
by Lipschitz continuity of $g$. The last bound does not depend on $\mathbf{z}%
\in \mathbf{R}^{d}$. Hence using Assumption 5(ii), we conclude 
\begin{equation*}
\sup_{\mathbf{z}\in \mathbf{R}^{d}}\left\vert \bar{g}_{n}(\mathbf{z})-\hat{g}%
_{n}(\mathbf{z})\right\vert =O_{P}\left( n^{-1/2}\varepsilon
_{n}^{-1}\right) ,
\end{equation*}%
where the convergence is uniform over $h\in H$. Therefore, as $M\rightarrow
\infty ,$%
\begin{equation*}
\lim_{n\rightarrow \infty }\sup_{h\in H}P_{n,h}\left\{ \sup_{c\in \lbrack
-M_{1},M_{1}]}\left\vert \bar{B}_{M_{1}}(c;a)-\hat{B}_{M_{1}}(c;a)\right%
\vert >Mn^{-1/2}\varepsilon _{n}^{-1}\right\} \rightarrow 0.
\end{equation*}

Since $g$ is Lipschitz, there exists $C>0$ such that for all $n\geq 1$, for
any $\mathbf{z},\mathbf{w}\in \mathbf{R}^{d}$,%
\begin{equation*}
\left\vert \bar{g}_{n}(\mathbf{z})-\bar{g}_{n}(\mathbf{w})\right\vert \leq
C||\mathbf{z-w}||.
\end{equation*}%
Hence by Assumptions 4(ii) and 5(i), we have%
\begin{equation*}
\lim_{n\rightarrow \infty }\sup_{h\in H}P_{n,h}\left\{ \sup_{c\in \lbrack
-M_{1},M_{1}]}\left\vert \tilde{B}_{M_{1}}(c;a)-\bar{B}_{M_{1}}(c;a)\right%
\vert >Mn^{-1/2}\right\} \rightarrow 0,
\end{equation*}%
as $M\rightarrow \infty $.

Now we show that as $M\rightarrow \infty $%
\begin{equation}
\lim_{n\rightarrow \infty }P\left\{ \sup_{c\in \lbrack
-M_{1},M_{1}]}\left\vert B_{M_{1}}^{\ast }(c;a)-\tilde{B}_{M_{1}}(c;a)\right%
\vert >M(L^{-1/2}+n^{-1/2})\right\} \rightarrow 0.  \label{cv4}
\end{equation}%
(Note that $P$ above denotes the joint distribution of the simulated
quantities $\{\xi _{i}\}_{i=1}^{L}$, and hence does not depend on $h\in H$.
Thus the convergence above is trivially uniform in $h\in H$.) First, define $%
f_{n}(\xi ;c,\mathbf{r})\equiv \tau _{M_{1}}\left( a|\bar{g}_{n}(\xi +%
\mathbf{r})-\bar{g}_{n}(\mathbf{r})+c|\right) $ and $\mathcal{J}_{n}\equiv
\{f_{n}(\cdot ;c,\mathbf{r}):(c,\mathbf{r})\in \lbrack -M_{1},M_{1}]\times
\lbrack -M_{1},M_{1}]^{d})\}$. The class $\mathcal{J}$ is uniformly bounded,
and $f(\xi ;c,\mathbf{r})$ is Lipschitz continuous in $(c,\mathbf{r})\in
\lbrack -M_{1},M_{1}]\times \lbrack -M_{1},M_{1}]^{d})$. Using the maximal
inequality (e.g. Theorems 2.14.2 (p.240) and 2.7.11 (p.164) of van der Vaart
and Wellner (1996)), we find that for some $C_{M_{1}}>0$ that depends only
on $M_{1}>0,$%
\begin{equation}
\mathbf{E}\left[ \sup_{c\in \lbrack -M_{1},M_{1}]}\left\vert B_{M_{1}}^{\ast
}(c;a)-\tilde{B}_{M_{1}}(c;a)\right\vert \right] \leq C_{M_{1}}\left\{
L^{-1/2}+n^{-1/2}\right\} .  \label{conv3}
\end{equation}%
Hence the convergence in (\ref{cv4}) follows. Thus the proof is complete. $%
\blacksquare $\bigskip

\noindent \textsc{Lemma A6}: \textit{Suppose that Assumptions} 1(i)\textit{\
and }4 \textit{hold. Then as }$n\rightarrow \infty ,$%
\begin{equation*}
\sup_{c\mathbf{\in }[-M_{1},M_{1}]}\left\vert B_{M_{1}}^{\ast
}(c;a)-B_{M_{1}}(c;a)\right\vert \rightarrow 0.
\end{equation*}

\noindent \textsc{Proof}: Since $g$ is Lipschitz continuous, the convergence 
\begin{equation*}
\bar{g}_{n}(\mathbf{z})\rightarrow \tilde{g}_{0}(\mathbf{z})\text{, as }%
n\rightarrow \infty ,
\end{equation*}%
is uniform over $\mathbf{z}$ in any given bounded subset of $\mathbf{R}^{d}$%
. (See Shapiro (1990), p.484.) Then 
\begin{eqnarray*}
&&\sup_{c\mathbf{\in }[-M_{1},M_{1}]}\left\vert B_{M_{1}}^{\ast
}(c;a)-B_{M_{1}}(c;a)\right\vert \\
&\leq &\sup_{c\mathbf{\in }[-M_{1},M_{1}]}\sup_{\mathbf{r}\in \lbrack
-M_{1},M_{1}]^{d}}\left\vert 
\begin{array}{c}
\mathbf{E}\left[ \tau _{M_{1}}\left( a\left\vert \bar{g}_{n}(Z+\mathbf{r})-%
\bar{g}_{n}(\mathbf{r})+c\right\vert \right) \right] \\ 
-\mathbf{E}\left[ \tau _{M_{1}}\left( a\left\vert \tilde{g}_{0}(Z+\mathbf{r}%
)-\tilde{g}_{0}(\mathbf{r})+c\right\vert \right) \right]%
\end{array}%
\right\vert \rightarrow 0,
\end{eqnarray*}%
as $n\rightarrow \infty $, because the domains of supremums above are
bounded in a finite dimensional space. $\blacksquare $\bigskip

\noindent \textsc{Lemma A7}: \textit{Suppose that Assumptions} 1(i), 4,%
\textit{\ and }5 \textit{hold. Then there exists }$M_{0}$\textit{\ such that
for any} $M_{1}>M_{0},$\textit{\ }$\varepsilon >0,$\textit{\ }$b>0$, \textit{%
and any }$a>0,$%
\begin{equation*}
\sup_{h\in H}P_{n,h}\left\{ \left\vert \hat{c}_{M_{1}}(a)-c_{M_{1}}^{\ast
}(a)\right\vert >\varepsilon \right\} \rightarrow 0,
\end{equation*}%
\textit{as} $n,L\rightarrow \infty $ \textit{jointly.}\bigskip

\noindent \textsc{Proof:} Let the Hausdorff distance between the two subsets 
$E_{1}$ and $E_{2}$ of $\mathbf{R}$ be denoted by $d_{H}(E_{1},E_{2})$.
First we show that 
\begin{equation}
d_{H}(E_{M_{1}}^{\ast }(a),\hat{E}_{M_{1}}(a))\rightarrow _{P}0,  \label{HDC}
\end{equation}%
as $n\rightarrow \infty $ and $L\rightarrow \infty $ uniformly over $h\in H$%
. For this, we use arguments in the proof of Theorem 3.1 of Chernozhukov,
Hong and Tamer (2007). Fix $\varepsilon \in (0,1)$ and let $E_{M_{1}}^{\ast
\varepsilon }(a)\equiv \{x\in \lbrack -M_{1},M_{1}]:\inf_{y\in
E_{M_{1}}^{\ast }}|x-y|\leq \varepsilon \}$. It suffices for (\ref{HDC}) to
show that for any $\varepsilon >0,$%
\begin{equation*}
\begin{tabular}{l}
(a) $\inf_{h\in H}P_{n,h}\left\{ \text{sup}_{c\in E_{M_{1}}^{\ast }(a)}\hat{B%
}_{M_{1}}(c;a)\leq \text{inf}_{c\in \lbrack -M_{1},M_{1}]}\hat{B}%
_{M_{1}}(c;a)+\eta _{n,L}\right\} \rightarrow 1,$ \\ 
(b) $\inf_{h\in H}P_{n,h}\left\{ \text{sup}_{c\in \hat{E}%
_{M_{1}}(a)}B_{M_{1}}^{\ast }(c;a)<\text{inf}_{c\in \lbrack
-M_{1},M_{1}]\backslash E_{M_{1}}^{\ast \varepsilon }(a)}B_{M_{1}}^{\ast
}(c;a)\right\} \rightarrow 1$,%
\end{tabular}%
\end{equation*}%
as $n,L\rightarrow \infty $ jointly. This is because (a) implies $\inf_{h\in
H}P_{n,h}\{E_{M_{1}}^{\ast }(a)\subset \hat{E}_{M_{1}}(a)\}\rightarrow 1$
and (b) implies that $\inf_{h\in H}P_{n,h}\{\hat{E}_{M_{1}}(a)\cap
([-M_{1},M_{1}]\backslash E_{M_{1}}^{\ast \varepsilon }(a))=\varnothing
\}\rightarrow 1$ so that $\inf_{h\in H}P_{n,h}\{\hat{E}_{M_{1}}(a)\subset
E_{M_{1}}^{\ast \varepsilon }(a)\}\rightarrow 1,$ and hence for any $%
\varepsilon >0$, 
\begin{equation*}
\text{sup}_{h\in H}P_{n,h}\left\{ d_{H}(E_{M_{1}}^{\ast }(a),\hat{E}%
_{M_{1}}(a))>\varepsilon \right\} \rightarrow 0,\ \text{as\ }n,L\rightarrow
\infty \ \text{jointly,}
\end{equation*}%
delivering (\ref{HDC}).

We focus on (a). Note that%
\begin{eqnarray*}
\text{sup}_{c\in E_{M_{1}}^{\ast }(a)}\hat{B}_{M_{1}}(c;a) &=&\text{sup}%
_{c\in E_{M_{1}}^{\ast }(a)}B_{M_{1}}^{\ast
}(c;a)+o_{P}(L^{-1/2}+n^{-1/2}\varepsilon _{n}^{-1}) \\
&\leq &\text{inf}_{c\in \lbrack -M_{1},M_{1}]}\hat{B}%
_{M_{1}}(c;a)+o_{P}(L^{-1/2}+n^{-1/2}\varepsilon _{n}^{-1}),
\end{eqnarray*}%
where the equality follows from Lemma A5, and the inequality follows by the
definition of $E_{M_{1}}^{\ast }(a)$. From this (a) follows because $\eta
_{n,L}\varepsilon _{n}\sqrt{n}\rightarrow \infty $ as $n\rightarrow \infty $
and $\eta _{n,L}\sqrt{L}\rightarrow \infty $ as $L\rightarrow \infty $.

Now let us turn to (b). Fix $\varepsilon >0$. Uniformly over $h\in H,$%
\begin{eqnarray}
\text{sup}_{c\in \hat{E}_{M_{1}}(a)}B_{M_{1}}^{\ast }(c;a) &\leq &\text{sup}%
_{c\in \hat{E}_{M_{1}}(a)}\hat{B}_{M_{1}}(c;a)+o_{P}(1)  \label{ineqs3} \\
&\leq &\inf_{c\in \lbrack -M_{1},M_{1}]}\hat{B}_{M_{1}}(c;a)+o_{P}(1)  \notag
\\
&\leq &\inf_{c\in \lbrack -M_{1},M_{1}]}B_{M_{1}}^{\ast }(c;a)+o_{P}(1), 
\notag
\end{eqnarray}%
where the second inequality follows by the definition of $\hat{E}_{M_{1}}(a)$
and the third inequality is due to $\eta _{n,L}\rightarrow 0$ as $%
n,L\rightarrow \infty $ and Lemma A5. By the definition of $E_{M_{1}}^{\ast
}(a)$, we have 
\begin{equation*}
0\leq \inf_{c\in \lbrack -M_{1},M_{1}]}B_{M_{1}}^{\ast }(c;a)<\inf_{c\in
\lbrack -M_{1},M_{1}]\backslash E_{M_{1}}^{\ast \varepsilon
}(a)}B_{M_{1}}^{\ast }(c;a).
\end{equation*}%
Hence we obtain (b). Thus we obtain (\ref{HDC}).

Now we show that as $n\rightarrow \infty ,$%
\begin{equation}
d_{H}(E_{M_{1}}^{\ast }(a),E_{M_{1}}(a))\rightarrow 0.  \label{HDC2}
\end{equation}%
Similarly as before, it suffices to note that 
\begin{eqnarray*}
\text{sup}_{c\in E_{M_{1}}(a)}B_{M_{1}}(c;a) &\leq &\text{inf}_{c\in \lbrack
-M_{1},M_{1}]}B_{M_{1}}^{\ast }(c;a)+o(1)\text{ and} \\
\text{sup}_{c\in E_{M_{1}}^{\ast }(a)}B_{M_{1}}(c;a) &<&\text{inf}_{c\in
\lbrack -M_{1},M_{1}]\backslash E_{M_{1}}^{\varepsilon }(a)}B_{M_{1}}(c;a).
\end{eqnarray*}%
The first inequality follows by the definition of $E_{M_{1}}(a)$ and Lemma
A6. The second inequality follows by Lemma A6 and the definition of $%
E_{M_{1}}^{\ast }(a)$ as in (\ref{ineqs3}). Thus we obtain (\ref{HDC2}). We
combine (\ref{HDC}) with (\ref{HDC2}) to conclude that%
\begin{equation}
d_{H}(E_{M_{1}}(a),\hat{E}_{M_{1}}(a))\rightarrow _{P}0.  \label{HDC3}
\end{equation}

For the main conclusion of the lemma, observe that $\left\vert \hat{c}%
_{M_{1}}(a)-c_{M_{1}}^{\ast }(a)\right\vert $ is equal to%
\begin{equation*}
\left\vert \sup \hat{E}_{M_{1}}(a)-\sup
E_{M_{1}}(a)\right\vert,
\end{equation*}%
which we can write as%
\begin{eqnarray*}
\left\vert \sup_{y\in \hat{E}_{M_{1}}(a)}\left\{ y-\sup
E_{M_{1}}(a)\right\} \right\vert 
=\left\vert \sup_{y\in \hat{E}_{M_{1}}(a)}\inf_{x\in
E_{M_{1}}(a)}\left( y-x\right)\right\vert.
\end{eqnarray*}%
We can interchange the supremum and the infimum using the fact that the sets $%
\hat{E}_{M_{1}}(a)$ and $E_{M_{1}}(a)$ are compact sets and using a version
of minimax theorem (e.g. Lemma A.3 of Puhalskii and Spokoiny (1998)). (Note
that the compactness of $\hat{E}_{M_{1}}(a)$ and $E_{M_{1}}(a)$ follows from
Assumption 4(ii).) Using the fact that $z=(z)_{+}-(z)_{-}$, where $%
(z)_{+}=\max (z,0)$ and $(z)_{-}=\max (-z,0)$, and applying the minimax
theorem, we bound the last term by%
\begin{eqnarray*}
	\sup_{y\in \hat{E}_{M_{1}}(a)}\inf_{x\in
		E_{M_{1}}(a)}\left( y-x\right) _{+} + \sup_{y\in E_{M_{1}}(a)}\inf_{x\in \hat{E}%
		_{M_{1}}(a)}\left( y-x\right)_{-}.
\end{eqnarray*}%
The sum above is bounded by $2d_{H}(E_{M_{1}}(a),\hat{E}_{M_{1}}(a))$. The
desired result follows from (\ref{HDC3}). $\blacksquare $\bigskip

\noindent \textsc{Proof of Theorem 2: }Fix $M>0$ and $\varepsilon >0$, and
take large $M_{1}\geq M$ such that%
\begin{eqnarray}
&&\sup_{\mathbf{r}\in \mathbf{R}^{d}}\mathbf{E}\left[ \tau _{M}(|\tilde{g}%
_{0}(Z+\mathbf{r})-\tilde{g}_{0}(\mathbf{r})+c_{M_{1}}^{\ast }(1)|)\right]
\label{ineq} \\
&\leq &\sup_{\mathbf{r}\in \lbrack -M_{1},M_{1}]^{d}}\mathbf{E}\left[ \tau
_{M}(|\tilde{g}_{0}(Z+\mathbf{r})-\tilde{g}_{0}(\mathbf{r})+c_{M_{1}}^{\ast
}(1)|)\right] +\varepsilon .  \notag
\end{eqnarray}%
This is possible for any choice of $\varepsilon >0$ because $\tau _{M}(\cdot
)$ and $\tilde{g}_{0}(\cdot )$ are Lipschitz continuous (recall Assumption
4(ii) and Lemma 1(ii)) and bounded by $M$. Note that%
\begin{eqnarray}
&&\sup_{h\in H_{n,b}}\mathbf{E}_{h}\left[ \tau _{M}(\sqrt{n}|\hat{\theta}-g(%
\mathbf{\beta }_{n}(h))|)\right]  \label{sup} \\
&=&\sup_{h\in H_{n,b}}\mathbf{E}_{h}\left[ \tau _{M}(\sqrt{n}|g(\mathbf{\hat{%
\beta}})+\hat{c}_{M_{1}}(1)/\sqrt{n}-g(\mathbf{\beta }_{n}(h))|)\right] 
\notag \\
&\leq &\sup_{h\in H_{n,b}}\mathbf{E}_{h}\left[ \tau _{M}(|g(\sqrt{n}\{%
\mathbf{\hat{\beta}}-\mathbf{\beta }_{n}(h)\}+\mathbf{r}_{n}(h))+\hat{c}%
_{M_{1}}(1)|)\right] ,  \notag
\end{eqnarray}%
where $\mathbf{r}_{n}(h)\equiv \sqrt{n}(\mathbf{\beta }_{n}(h)-g(\mathbf{%
\beta }_{n}(h)))$. Note that for each $h\in H_{n,b}$,%
\begin{eqnarray*}
\mathbf{r}_{n}(h) &=&\sqrt{n}\{\mathbf{\beta }_{n}(h)-\mathbf{\beta }_{n}(0)+%
\mathbf{\beta }_{n}(0)\}-\sqrt{n}g(\mathbf{\beta }_{n}(h)-\mathbf{\beta }%
_{n}(0)+\mathbf{\beta }_{n}(0)) \\
&=&\{\sqrt{n}\mathbf{\beta }_{0}+\mathbf{\tilde{r}}_{n}(h)-g(\sqrt{n}\mathbf{%
\beta }_{0}+\mathbf{\tilde{r}}_{n}(h))\},
\end{eqnarray*}%
where $\mathbf{\tilde{r}}_{n}(h)\equiv \sqrt{n}\{\mathbf{\beta }_{n}(h)-%
\mathbf{\beta }_{n}(0)\}$ and $\sup_{h\in H_{n,b}}||\mathbf{\tilde{r}}%
_{n}(h)||\leq b||\mathbf{\dot{\beta}}^{\ast }||+o(1)$ by the definition of $%
h\in H_{n,b}$. Using Assumption 5, and using the fact that $Z$ is a
continuous random vector, we find that%
\begin{equation*}
\sup_{t\in \mathbf{R}^{d}}\sup_{\mathbf{r}\in \mathbf{R}^{d}}\sup_{h\in
H}\left\vert 
\begin{array}{c}
P_{n,h}\left\{ \sqrt{n}\{\mathbf{\hat{\beta}}-\mathbf{\beta }_{n}(h)\}+%
\mathbf{r}+\hat{c}_{M_{1}}(1)\leq t\right\} \\ 
-P\left\{ Z+\mathbf{r}+c_{M_{1}}^{\ast }(1)\leq t\right\}%
\end{array}%
\right\vert \rightarrow 0,
\end{equation*}%
as $n\rightarrow \infty $. Therefore,%
\begin{equation*}
\sup_{h\in H_{n,b}}\left\vert 
\begin{array}{c}
\mathbf{E}_{h}\left[ \tau _{M}(|g(\sqrt{n}\{\mathbf{\hat{\beta}}-\mathbf{%
\beta }_{n}(h)\}+\mathbf{r}_{n}(h))+\hat{c}_{M_{1}}(1)|)\right] \\ 
-\mathbf{E}\left[ \tau _{M}(|g(Z+\mathbf{r}_{n}(h))+c_{M_{1}}^{\ast }(1)|)%
\right]%
\end{array}%
\right\vert \rightarrow 0,
\end{equation*}%
as $n\rightarrow \infty $. Let $A_{M}\equiv \tau ^{-1}\left( [0,M]\right) $
which is bounded in $[0,\infty )$ by Assumption 4(i). We take $M_{2}\geq
M_{1}$ and write 
\begin{eqnarray*}
&&\mathbf{E}\left[ \tau _{M}(|g(Z+\mathbf{r}_{n}(h))+c_{M_{1}}^{\ast }(1)|)%
\right] \\
&\leq &\mathbf{E}\left[ \tau _{M}(|g(Z+\mathbf{r}_{n}(h))+c_{M_{1}}^{\ast
}(1)|)1\left\{ ||Z||\leq M_{2}\right\} \right] +MP\left\{
||Z||>M_{2}\right\} ,
\end{eqnarray*}%
where the leading expectation in the second line can be rewritten as%
\begin{equation}
\mathbf{E}\left[ \tau _{M}\left( \left\vert 
\begin{array}{c}
g\left( Z+\sqrt{n}\mathbf{\beta }_{0}+\mathbf{\tilde{r}}_{n}(h)\right) \\ 
-g(\sqrt{n}\mathbf{\beta }_{0}+\mathbf{\tilde{r}}_{n}(h))+c_{M_{1}}^{\ast
}(1)%
\end{array}%
\right\vert \right) 1\left\{ ||Z||\leq M_{2}\right\} \right] .  \label{exp}
\end{equation}%
Since $g$ is Lipschitz, $\sup_{h\in H_{n,b}}||\mathbf{\tilde{r}}%
_{n}(h)||\leq b||\mathbf{\dot{\beta}}^{\ast }||+o(1)$, and the convergence
of 
\begin{equation*}
g(\mathbf{z}+\sqrt{n}\mathbf{\beta }_{0})-g(\sqrt{n}\mathbf{\beta }%
_{0})\rightarrow \tilde{g}_{0}(\mathbf{z})
\end{equation*}%
is uniform over $\mathbf{z}$ in any bounded set by Lemma 1(iii), we find
that the expectation in (\ref{exp}) converges to%
\begin{equation*}
\mathbf{E}\left[ \tau _{M}(|\tilde{g}_{0}(Z+\mathbf{\dot{\beta}}(h))-\tilde{g%
}_{0}(\mathbf{\dot{\beta}}(h))+c_{M_{1}}^{\ast }(1))|1\left\{ ||Z||\leq
M_{2}\right\} \right] ,
\end{equation*}%
uniformly in $h\in H_{n,b}$ as $n\rightarrow \infty $. Thus, we conclude
that 
\begin{eqnarray*}
&&\underset{n\rightarrow \infty }{\text{limsup}}\sup_{h\in H_{n,b}}\mathbf{E}%
_{h}\left[ \tau _{M}(\sqrt{n}|\hat{\theta}-g(\mathbf{\beta }_{n}(h))|)\right]
\\
&\leq &\sup_{\mathbf{r}\in \mathbf{R}^{d}}\mathbf{E}\left[ \tau _{M}(|\tilde{%
g}_{0}(Z+\mathbf{r})-\tilde{g}_{0}(\mathbf{r})+c_{M_{1}}^{\ast
}(1))|1\left\{ ||Z||\leq M_{2}\right\} \right] +MP\left\{
||Z||>M_{2}\right\} .
\end{eqnarray*}%
As we send $M_{2}\uparrow \infty $, the last sum vanishes and the leading
supremum becomes%
\begin{eqnarray*}
&&\sup_{\mathbf{r}\in \mathbf{R}^{d}}\mathbf{E}\left[ \tau _{M}(|\tilde{g}%
_{0}(Z+\mathbf{r})-\tilde{g}_{0}(\mathbf{r})+c_{M_{1}}^{\ast }(1))|\right] \\
&\leq &\sup_{\mathbf{r}\in \lbrack -M_{1},M_{1}]^{d}}\mathbf{E}\left[ \tau
_{M}(|\tilde{g}_{0}(Z+\mathbf{r})-\tilde{g}_{0}(\mathbf{r})+c_{M_{1}}^{\ast
}(1)|)\right] +\varepsilon ,
\end{eqnarray*}%
by (\ref{ineq}). Since $M_{1}\geq M$, the last supremum is bounded by%
\begin{eqnarray*}
&&\sup_{\mathbf{r}\in \lbrack -M_{1},M_{1}]^{d}}\mathbf{E}\left[ \tau
_{M_{1}}(|\tilde{g}_{0}(Z+\mathbf{r})-\tilde{g}_{0}(\mathbf{r}%
)+c_{M_{1}}^{\ast }(1)|)\right] \\
&=&\inf_{-M_{1}\leq c\leq M_{1}}\sup_{\mathbf{r}\in \lbrack
-M_{1},M_{1}]^{d}}\mathbf{E}\left[ \tau _{M_{1}}(|\tilde{g}_{0}(Z+\mathbf{r}%
)-\tilde{g}_{0}(\mathbf{r})+c|)\right] ,
\end{eqnarray*}%
where the equality follows by the definition of $c_{M_{1}}^{\ast }(1)$.
Since the choice of $\varepsilon $ and $M_{1}$ was arbitrary and $\mathbf{E}%
\left[ \tau _{M_{1}}(|\tilde{g}_{0}(Z+\mathbf{r})-\tilde{g}_{0}(\mathbf{r}%
)+c|)\right] $ is uniformly continuous in $\mathbf{r}\in \mathbf{R}^{d},$
sending $M_{1}\uparrow \infty $ (along with $\varepsilon \downarrow 0$), and
then sending $M\uparrow \infty $, we obtain the desired result. $%
\blacksquare $\bigskip

\noindent \textsc{Proof of Theorem 3:} As in the proof of Lemma 3, we choose 
$\mathbf{r}\in \mathbf{R}^{d}$ so that for some $h^{\prime }\in H$, $\mathbf{%
r=\dot{\beta}(}h^{\prime }\mathbf{)}$. Fix $b/2\geq ||h^{\prime }||\cdot ||%
\mathbf{\dot{\beta}}^{\ast }||$. Define%
\begin{eqnarray*}
H_{n,b,1}^{\ast } &\equiv &\{h\in H_{n,b}^{\ast }:g(\mathbf{\beta }%
_{n}(h+h^{\prime }))\geq \bar{x}\},\text{ and} \\
H_{n,b,2}^{\ast } &\equiv &\{h\in H_{n,b}^{\ast }:g(\mathbf{\beta }%
_{n}(h+h^{\prime }))\leq \bar{x}\},
\end{eqnarray*}%
where we recall $H_{n,b}^{\ast }\equiv \{h\in H_{n,b}:\langle h,h^{\prime
}\rangle =0\}$.

First, suppose that $g(\mathbf{\beta }_{0})>\bar{x}$. Note that 
\begin{eqnarray*}
&&\underset{n\rightarrow \infty }{\text{liminf}}\sup_{h\in H_{n,b}}\mathbf{E}%
_{h}\left[ \tau (|\sqrt{n}\{\hat{\theta}-f(g(\mathbf{\beta }_{n}(h)))\}|)%
\right] \\
&\geq &\ \underset{n\rightarrow \infty }{\text{liminf}}\sup_{h\in
H_{n,b/2,1}^{\ast }}\mathbf{E}_{h+h^{\prime }}\left[ \tau _{M}\left( |\sqrt{n%
}\{\hat{\theta}+a_{1}\bar{x}-a_{1}g(\mathbf{\beta }_{n}(h+h^{\prime
}))\}|\right) \right] \\
&\geq &\ \underset{n\rightarrow \infty }{\text{liminf}}\sup_{h\in
H_{n,b/2,1}^{\ast }}\mathbf{E}_{h+h^{\prime }}\left[ \tau _{M}\left( |\sqrt{n%
}\{\tilde{\theta}_{1}-a_{1}g(\mathbf{\beta }_{n}(h+h^{\prime }))\}|\right) %
\right] ,
\end{eqnarray*}%
where $\tilde{\theta}_{1}\equiv \hat{\theta}+a_{1}\bar{x}$. Let $\tilde{V}%
_{n,h,1}\equiv \sqrt{n}\{\tilde{\theta}_{1}-a_{1}g(\mathbf{\beta }_{n}(h))\}$%
, $h(\mathbf{p})$, $\mathbf{p}=(p_{i})_{i=1}^{m}\in \mathbf{R}^{m}$, and $%
F_{\lambda }(\mathbf{p})$ be as in the proof of Lemma 3, so that we have%
\begin{eqnarray*}
&&\underset{n\rightarrow \infty }{\text{liminf}}\sup_{h\in H_{n,b/2,1}^{\ast
}}\mathbf{E}_{h+h^{\prime }}\left[ \tau _{M}\left( |\sqrt{n}\{\tilde{\theta}%
_{1}-a_{1}g(\mathbf{\beta }_{n}(h+h^{\prime }))|\right) \right] \\
&\geq &\ \int \underset{n\rightarrow \infty }{\text{liminf}}\ \mathbf{E}_{h(%
\mathbf{p})+h^{\prime }}\left[ \tau _{M}\left( |\tilde{V}_{n,h(\mathbf{p}%
)+h^{\prime },1}|\right) \right] 1\left\{ h(\mathbf{p})\in H_{n,b,1}\right\}
dF_{\lambda }(\mathbf{p}).
\end{eqnarray*}%
Let $R_{n}(h)\equiv g(\mathbf{\beta }_{n}(h+h^{\prime }))-g(\mathbf{\beta }%
_{0})$, and observe that for all $h\in H,$%
\begin{eqnarray*}
R_{n}(h) &=&g(\mathbf{\beta }_{n}(h+h^{\prime })-\mathbf{\beta }%
_{n}(h^{\prime })+\mathbf{\beta }_{n}(h^{\prime }))-g(\mathbf{\beta }%
_{n}(h^{\prime })) \\
&&+g(\mathbf{\beta }_{n}(h^{\prime }))-g(\mathbf{\beta }_{0}) \\
&=&g((\mathbf{\dot{\beta}}(h+h^{\prime }))/\sqrt{n}+\mathbf{\beta }%
_{n}(h^{\prime }))-g(\mathbf{\beta }_{n}(h^{\prime })) \\
&&+g(\mathbf{\beta }_{n}(h^{\prime }))-g(\mathbf{\beta }_{0})+o(1/\sqrt{n}),
\end{eqnarray*}%
as $n\rightarrow \infty .$ Since the map $g$ is Lipshitz continuous and $%
\mathbf{\beta }_{n}(h^{\prime })=\mathbf{\beta }_{0}+O(1/\sqrt{n})$), we
deduce that for each $h\in H$,%
\begin{equation}
\left\vert R_{n}(h)\right\vert \rightarrow 0,  \label{cv2}
\end{equation}%
as $n\rightarrow \infty $. This means that given that $g(\mathbf{\beta }%
_{0})>\bar{x}$, we have 
\begin{equation*}
1\left\{ h(\mathbf{p})\in H_{n,b,1}^{\ast }\right\} \rightarrow 1\left\{ h(%
\mathbf{p})\in H_{b}^{\ast }\right\} ,
\end{equation*}%
as $n\rightarrow \infty $.

Following the same arguments as in the proofs of Lemma 3 and Theorem 1, we
deduce that%
\begin{eqnarray*}
&&\sup_{b\in \lbrack 0,\infty )}\underset{n\rightarrow \infty }{\text{liminf}%
}\sup_{h\in H_{n,b,1}^{\ast }}\mathbf{E}_{h}\left[ \tau _{M}\left( |\sqrt{n}%
\{\tilde{\theta}_{1}-a_{1}g(\mathbf{\beta }_{n}(h))|\right) \right] \\
&\geq &\inf_{c\in \mathbf{R}}\sup_{\mathbf{r}\in \mathbf{R}^{d}}\mathbf{E}%
\left[ \tau _{M}(|a_{1}||\tilde{g}_{0}(Z+\mathbf{r})-\tilde{g}_{0}(\mathbf{r}%
)+c|)\right] -\varepsilon .
\end{eqnarray*}

Second, suppose that $g(\mathbf{\beta }_{0})<\bar{x}$. Using similar
arguments, we obtain the result that%
\begin{eqnarray*}
&&\sup_{b\in \lbrack 0,\infty )}\underset{n\rightarrow \infty }{\text{liminf}%
}\sup_{h\in H_{n,b,2}^{\ast }}\mathbf{E}_{h}\left[ \tau _{M}\left( |\sqrt{n}%
\{\tilde{\theta}_{2}-a_{2}g(\mathbf{\beta }_{n}(h))|\right) \right] \\
&\geq &\inf_{c\in \mathbf{R}}\sup_{\mathbf{r}\in \mathbf{R}^{d}}\mathbf{E}%
\left[ \tau _{M}(|a_{2}||\tilde{g}_{0}(Z+\mathbf{r})-\tilde{g}_{0}(\mathbf{r}%
)+c|)\right] -\varepsilon ,
\end{eqnarray*}%
where $\tilde{\theta}_{1}\equiv \hat{\theta}+a_{2}\bar{x}$.

Finally, assume that $g(\mathbf{\beta }_{0})=\bar{x}$. Then 
\begin{eqnarray*}
H_{n,b,1}^{\ast } &=&\{h\in H_{n,b}^{\ast }:R_{n}(h)\geq 0\},\text{ and} \\
H_{n,b,2}^{\ast } &=&\{h\in H_{n,b}^{\ast }:R_{n}(h)\leq 0\}.
\end{eqnarray*}%
By (\ref{cv2}), we have for each $h\in H$, as $n\rightarrow \infty $,%
\begin{eqnarray}
1\left\{ h\in H_{n,b,1}^{\ast }\right\} &\rightarrow &1\left\{ h\in
H_{b}\right\} \text{ and}  \label{cv3} \\
1\left\{ h\in H_{n,b,2}^{\ast }\right\} &\rightarrow &1\left\{ h\in
H_{b}\right\} ,  \notag
\end{eqnarray}%
where $H_{b}\equiv \{h\in H:||\mathbf{\dot{\beta}}(h)||\leq b\}$. Note that 
\begin{equation*}
\sup_{h\in H_{n,b}}\mathbf{E}_{h}\left[ \tau _{M}\left( |\sqrt{n}\{\hat{%
\theta}-f(g(\mathbf{\beta }_{n}(h)))|\right) \right] \geq
\max_{l=1,2}\sup_{h\in H_{n,b/2,l}^{\ast }}\mathbf{E}_{h}\left[ \tau
_{M}\left( |\sqrt{n}\{\tilde{\theta}_{l}-a_{l}g(\mathbf{\beta }%
_{n}(h))|\right) \right] .
\end{equation*}%
Using (\ref{cv3}) and following the same arguments as before, we conclude
that%
\begin{eqnarray*}
&&\sup_{b\in \lbrack 0,\infty )}\underset{n\rightarrow \infty }{\text{liminf}%
}\sup_{h\in H_{n,b}}\mathbf{E}_{h}[\tau (|\sqrt{n}\{\hat{\theta}-f(g(\mathbf{%
\beta }_{n}(h)))\}|)] \\
&\geq &\max_{l=1,2}\inf_{c\in \mathbf{R}}\sup_{\mathbf{r}\in \mathbf{R}^{d}}%
\mathbf{E}\left[ \tau _{M}(|a_{l}||\tilde{g}_{0}(Z+\mathbf{r})-\tilde{g}_{0}(%
\mathbf{r})+c|)\right] \\
&=&\inf_{c\in \mathbf{R}}\sup_{\mathbf{r}\in \mathbf{R}^{d}}\mathbf{E}\left[
\tau _{M}(\max \left\{ |a_{1}|,|a_{2}|\right\} |\tilde{g}_{0}(Z+\mathbf{r})-%
\tilde{g}_{0}(\mathbf{r})+c|)\right] ,
\end{eqnarray*}%
where the last equality follows because $\tau _{M}$ is an increasing
function. By sending $M\uparrow \infty $, we obtain the desired result. $%
\blacksquare $\bigskip

\noindent \textsc{Lemma A8:}\textit{\ Suppose that Assumptions 1(i) and 5
hold. Then,}%
\begin{equation*}
\inf_{h\in H}P_{n,h}\left\{ \hat{s}=s\right\} \rightarrow 1.
\end{equation*}%
\noindent \textsc{Proof:}\textit{\ }The proof can be straightforwardly
proceeded as the proof of Lemma A5 by dividing the proof into cases with $g(%
\mathbf{\beta }_{0})>\bar{x}$, $g(\mathbf{\beta }_{0})<\bar{x}$, and $g(%
\mathbf{\beta }_{0})=\bar{x}$, and applying Assumption A5(ii). The details
are omitted. $\blacksquare $\bigskip

\noindent \textsc{Proof of Theorem 4:} For any $M>0$, and $b\in \lbrack
0,\infty ),$ 
\begin{eqnarray*}
&&\underset{n\rightarrow \infty }{\text{limsup}}\sup_{h\in H_{n,b}}\mathbf{E}%
_{h}\left[ \tau _{M}(|\sqrt{n}\{\tilde{\theta}_{mx}-f(g(\mathbf{\beta }%
_{n}(h)))\}|)\right] \\
&=&\ \underset{n\rightarrow \infty }{\text{limsup}}\sup_{h\in H_{n,b}}%
\mathbf{E}_{h}\left[ \tau _{M}(|\sqrt{n}\{\tilde{\theta}_{mx}-f(g(\mathbf{%
\beta }_{n}(h)))\}|)1\left\{ \hat{s}=s\right\} \right] ,
\end{eqnarray*}%
by Lemma A8. We focus on the last limsup.

First, suppose that $g(\mathbf{\beta }_{0})>\bar{x}$. Then there exists $%
\varepsilon >0$, such that $g(\mathbf{\beta }_{0})>\bar{x}+\varepsilon $.
Since have for all $h\in H_{n,b}$,%
\begin{equation*}
||\mathbf{\beta }_{n}(h)-\mathbf{\beta }_{0}||\leq b/\sqrt{n}\text{,}
\end{equation*}%
we conclude that from some large $n$ on, for all $h\in H_{n,b}$, we have 
\begin{equation*}
g(\mathbf{\beta }_{n}(h))\geq \bar{x}.
\end{equation*}%
Hence 
\begin{eqnarray*}
&&\ \sup_{b\in \lbrack 0,\infty )}\underset{n\rightarrow \infty }{\text{%
limsup}}\sup_{h\in H_{n,b}}\mathbf{E}_{h}\left[ \tau _{M}(|\sqrt{n}\{\tilde{%
\theta}_{mx}-f(g(\mathbf{\beta }_{n}(h)))\}|)1\left\{ \hat{s}=s\right\} %
\right] \\
&=&\ \sup_{b\in \lbrack 0,\infty )}\underset{n\rightarrow \infty }{\text{%
limsup}}\sup_{h\in H_{n,b}}\mathbf{E}_{h}\left[ \tau _{M}(|a_{1}||\sqrt{n}\{%
\tilde{\theta}_{mx,1}-g(\mathbf{\beta }_{n}(h)))\}|)1\left\{ \hat{s}%
=s\right\} \right] \\
&=&\ \sup_{b\in \lbrack 0,\infty )}\underset{n\rightarrow \infty }{\text{%
limsup}}\sup_{h\in H_{n,b}}\mathbf{E}_{h}\left[ \tau _{M}(|a_{1}||\sqrt{n}\{%
\tilde{\theta}_{mx,1}-g(\mathbf{\beta }_{n}(h)))\}|)\right] ,
\end{eqnarray*}%
where $\tilde{\theta}_{mx,1}=\tilde{\theta}_{mx}/a_{1}$. By Assumption 5, we
have%
\begin{eqnarray*}
&&\sup_{h\in H_{n,b}}\mathbf{E}_{h}\left[ \tau _{M}(|a_{1}||\sqrt{n}\{\tilde{%
\theta}_{mx,1}-g(\mathbf{\beta }_{n}(h)))\}|)\right] \\
&\leq &\sup_{h\in H_{n,b}}\mathbf{E}_{h}\left[ \tau _{M}(\sqrt{n}|a_{1}||g(%
\mathbf{\hat{\beta}})+\hat{c}_{M_{1}}(s)/\sqrt{n}-g(\mathbf{\beta }_{n}(h))|)%
\right] .
\end{eqnarray*}%
Fix $\varepsilon >0$, choose $M_{1}\geq M$, and follow the proof of Theorem
2 to find that the limsup$_{n\rightarrow \infty }$ of the last supremum is
bounded by 
\begin{equation*}
\sup_{\mathbf{r}\in \lbrack -M_{1},M_{1}]^{d}}\mathbf{E}\left[ \tau
_{M_{1}}(|a_{1}||\tilde{g}_{0}(Z+\mathbf{r})-\tilde{g}_{0}(\mathbf{r}%
)+c_{M_{1}}^{\ast }(s)|)\right] +\varepsilon .
\end{equation*}%
By the definition of $c_{M_{1}}^{\ast }(s),$ the last supremum is equal to%
\begin{eqnarray*}
&&\inf_{c\in \lbrack -M_{1},M_{1}]}\sup_{\mathbf{r}\in \lbrack
-M_{1},M_{1}]^{d}}\mathbf{E}\left[ \tau _{M_{1}}(|a_{1}||\tilde{g}_{0}(Z+%
\mathbf{r})-\tilde{g}_{0}(\mathbf{r})+c|)\right] \\
&\leq &\inf_{c\in \lbrack -M_{1},M_{1}]}\sup_{\mathbf{r}\in \mathbf{R}^{d}}%
\mathbf{E}\left[ \tau (|a_{1}||\tilde{g}_{0}(Z+\mathbf{r})-\tilde{g}_{0}(%
\mathbf{r})+c|)\right] .
\end{eqnarray*}%
Sending $M_{1}\uparrow \infty $, we conclude that%
\begin{eqnarray*}
&&\sup_{b\in \lbrack 0,\infty )}\underset{n\rightarrow \infty }{\text{limsup}%
}\sup_{h\in H_{n,b}}\mathbf{E}_{h}\left[ \tau _{M}(|\sqrt{n}\{\tilde{\theta}%
_{mx}-f(g(\mathbf{\beta }_{n}(h)))\}|)1\left\{ \hat{s}=s\right\} \right] \\
&\leq &\inf_{c\in \mathbf{R}}\sup_{\mathbf{r}\in \mathbf{R}^{d}}\mathbf{E}%
\left[ \tau (|a_{1}||\tilde{g}_{0}(Z+\mathbf{r})-\tilde{g}_{0}(\mathbf{r}%
)+c|)\right] .
\end{eqnarray*}

Second, suppose that $g(\mathbf{\beta }_{0})<\bar{x}$. Then we can use the
same arguments as before to show the following:%
\begin{eqnarray*}
&&\sup_{b\in \lbrack 0,\infty )}\underset{n\rightarrow \infty }{\text{limsup}%
}\sup_{h\in H_{n,b}}\mathbf{E}_{h}\left[ \tau _{M}(|\sqrt{n}\{\tilde{\theta}%
_{mx}-f(g(\mathbf{\beta }_{n}(h)))\}|)1\left\{ \hat{s}=s\right\} \right] \\
&\leq &\inf_{c\in \mathbf{R}}\sup_{\mathbf{r}\in \mathbf{R}^{d}}\mathbf{E}%
\left[ \tau (|a_{2}||\tilde{g}_{0}(Z+\mathbf{r})-\tilde{g}_{0}(\mathbf{r}%
)+c|)\right] .
\end{eqnarray*}

Finally, suppose that $g(\mathbf{\beta }_{0})=\bar{x}$. Then note that $%
f(\cdot )/s$ with $s=\max \{|a_{1}|,|a_{2}|\}$ is a contraction mapping.
Hence%
\begin{eqnarray*}
&&\underset{n\rightarrow \infty }{\text{limsup}}\sup_{h\in H_{n,b}}\mathbf{E}%
_{h}\left[ \tau _{M}(|\sqrt{n}\{\tilde{\theta}_{mx}-f(g(\mathbf{\beta }%
_{n}(h)))\}|)1\left\{ \hat{s}=s\right\} \right] \\
&=&\ \underset{n\rightarrow \infty }{\text{limsup}}\sup_{h\in H_{n,b}}%
\mathbf{E}_{h}\left[ \tau _{M}(\sqrt{n}|f(g(\mathbf{\hat{\beta}})+\hat{c}%
_{M_{1}}(s)/\sqrt{n})-f(g(\mathbf{\beta }_{n}(h)))|)1\left\{ \hat{s}%
=s\right\} \right] \\
&\leq &\ \underset{n\rightarrow \infty }{\text{limsup}}\sup_{h\in H_{n,b}}%
\mathbf{E}_{h}\left[ \tau _{M}(\sqrt{n}s|g(\mathbf{\hat{\beta}})+\hat{c}%
_{M_{1}}(s)/\sqrt{n}-g(\mathbf{\beta }_{n}(h))|)1\left\{ \hat{s}=s\right\} %
\right] \\
&=&\ \underset{n\rightarrow \infty }{\text{limsup}}\sup_{h\in H_{n,b}}%
\mathbf{E}_{h}\left[ \tau _{M}(\sqrt{n}s|g(\mathbf{\hat{\beta}})+\hat{c}%
_{M_{1}}(s)/\sqrt{n}-g(\mathbf{\beta }_{n}(h))|)\right] ,
\end{eqnarray*}%
by Lemma A8, where the inequality above is due to $f(\cdot )/s$ being a
contraction mapping. We fix $\varepsilon >0$ and choose $M_{1}\geq M$ and
follow the proof of Theorem 2 to find that the sup$_{b\in \lbrack 0,\infty
)} $ of the last limsup is bounded by 
\begin{equation*}
\inf_{c\in \lbrack -M_{1},M_{1}]}\sup_{\mathbf{r}\in \mathbf{R}^{d}}\mathbf{E%
}\left[ \tau \left( s\left\vert \tilde{g}_{0}(Z+\mathbf{r})-\tilde{g}_{0}(%
\mathbf{r})+c\right\vert \right) \right] .
\end{equation*}%
By sending $M_{1}\uparrow \infty $, we obtain the desired bound. $%
\blacksquare $

\pagebreak

	\begin{center}
	\Large \textsc{Corrigendum to ``Local Asymptotic Minimax
		Estimation of Nonregular Parameters with Translation-Scale Equivariant Maps": [J. Multivariate Anal. 125 (2014) 136–158]}
	\bigskip
\end{center}

\vspace*{1ex minus 1ex}
\begin{center}
	Kyungchul Song\\
	\textit{Vancouver School of Economics, University of British Columbia}\\
	\bigskip
	
\end{center}

First, the proof of Theorem 1 contains a gap in the equation on page 151:
\begin{equation*}
	\min_{u\in \mathcal{T}_{K,N}}\max_{\mathbf{r}\in \mathcal{J}_{K}}\int \tilde{%
		g}_{0}(z+u)d\Lambda _{\mathbf{r}}(z)=\min_{u\in \mathcal{T}_{K,N}}\max_{%
		\mathbf{r}\in \mathcal{J}_{K}}\mathbf{E}\left[ \tau _{M}(|\tilde{g}_{0}(Z+%
	\mathbf{r})-\tilde{g}_{0}(\mathbf{r})+u|)\right].
\end{equation*}
(I thank Yoshiyasu Rai for pointing it out to me.) Theorem 1 still holds if we focus on convex loss functions, replacing Assumption 4 (i) on page 140 by the following:\medskip

\noindent \textbf{Assumption 4} (i) $\tau(\cdot)$ is increasing and convex on $[0,\infty)$, $\tau(0) = 0$, and there exists $\bar \tau$ such that $\tau^{-1}([0,y])$ is bounded in $[0,\infty)$ for all $0 < y < \bar \tau$.\medskip

Then Theorem 1 follows from Lemma 3 by Jensen's inequality, because 
\begin{eqnarray*}
	&& \inf_{F\in \mathcal{F}}\sup_{\mathbf{r}\in \mathbf{R}^{d}}\int 
	\mathbf{E}\left[ \tau (|\tilde{g}_{0}(Z+\mathbf{r})-\tilde{g}_{0}(\mathbf{r}%
	)+w|)\right] dF(w)\\
	&\ge& \inf_{F\in \mathcal{F}}\sup_{\mathbf{r}\in \mathbf{R}^{d}}
	\mathbf{E}\left[ \tau \left(\left|\tilde{g}_{0}(Z+\mathbf{r})-\tilde{g}_{0}(\mathbf{r}%
	)+\int w dF(w)\right|\right)\right] = \inf_{c \in \mathbf{R}} B(c;1).
\end{eqnarray*}
\medskip

Second, the last equality on page 155 in the proof of Theorem 2 as follows:
\begin{eqnarray*}
	&&\sup_{\mathbf{r}\in \lbrack -M_{1},M_{1}]^{d}}\mathbf{E}\left[ \tau
	_{M_{1}}(|\tilde{g}_{0}(Z+\mathbf{r})-\tilde{g}_{0}(\mathbf{r}%
	)+c_{M_{1}}^{\ast }(1)|)\right] \\
	&=&\inf_{-M_{1}\leq c\leq M_{1}}\sup_{\mathbf{r}\in \lbrack
		-M_{1},M_{1}]^{d}}\mathbf{E}\left[ \tau _{M_{1}}(|\tilde{g}_{0}(Z+\mathbf{r}%
	)-\tilde{g}_{0}(\mathbf{r})+c|)\right]
\end{eqnarray*}
assumes that the set $E_{M_1}(a)$ is convex, which is not guaranteed. (I thank Zheng Fang for pointing it out to me.) Note that the set $E_{M_1}(a)$ (defined in the first display on page 152) is compact. Hence the results of the paper including Theorem 2 follow once we redefine
\begin{eqnarray*}
	\hat c_{M_1}(a) &\equiv& \sup \hat E_{M_1}(a), \text{ in (3.4) on page 141, and }\\
	c_{M_1}^*(a) &\equiv& \sup E_{M_{1}}(a), \text{ in Line 3 on page 152.}
\end{eqnarray*}
Modifying the definitions using infimum in place of supremum works as well.\medskip

I apologize for any inconvenience caused by these gaps.\medskip

\end{document}